\documentclass[aos,preprint]{imsart}

\usepackage{xcolor}

\usepackage{hyperref}
\usepackage{indentfirst}
\usepackage{subcaption}
\usepackage{bm}

\newcommand{\argmin}{\operatorname{argmin}}

\usepackage{amsthm,amsmath,amssymb}
\usepackage[authoryear]{natbib}
\usepackage{caption, subcaption}
\usepackage{graphicx} 
\usepackage{algorithm}
\usepackage{comment}
\usepackage[noend]{algpseudocode}

\newcommand{\E}{\mathbb{E}}
\newcommand{\U}{\mathbb{U}}
\renewcommand{\P}{\mathbb{P}}

\newcommand{\bitem}{\begin{itemize}}
\newcommand{\eitem}{\end{itemize}}
\newcommand{\benum}{\begin{enumerate}}
\newcommand{\eenum}{\end{enumerate}}
\newcommand{\beq}{\begin{equation}}
\newcommand{\eeq}{\end{equation}}
\newcommand{\beqs}{\begin{equation*}}
\newcommand{\eeqs}{\end{equation*}}
\newcommand{\beas}{\begin{eqnarray*}}
\newcommand{\eeas}{\end{eqnarray*}}
\newcommand{\bea}{\begin{eqnarray}}
\newcommand{\eea}{\end{eqnarray}}
\newcommand{\bei}{\begin{itemize}}
\newcommand{\eei}{\end{itemize}}
\newcommand{\ben}{\begin{enumerate}}
\newcommand{\een}{\end{enumerate}}
\newcommand{\bet}{\begin{theorem}}
\newcommand{\eet}{\end{theorem}}
\newcommand{\bel}{\begin{lemma}}
\newcommand{\eel}{\end{lemma}}
\newcommand{\bep}{\begin{proposition}}
\newcommand{\eep}{\end{proposition}}
\newcommand{\bed}{\begin{definition}}
\newcommand{\eed}{\end{definition}}
\newcommand{\bec}{\begin{corollary}}
\newcommand{\eec}{\end{corollary}}
\newcommand{\bex}{\begin{example}}
\newcommand{\eex}{\end{example}}

\def\balign#1\ealign{\begin{align}#1\end{align}}
\def\baligns#1\ealigns{\begin{align*}#1\end{align*}}
\def\reals{\mathbb{R}} 
\def\integers{\mathbb{Z}} 
\def\naturals{\mathbb{N}} 

\newcommand{\ind}[1]{\mathbb{I}_{\{#1\}}}

\newcommand{\iid}{\stackrel{\mathrm{iid}}{\sim}}

\newtheorem{theorem}{Theorem}
\newtheorem{lemma}{Lemma}

\newtheorem{corollary}[theorem]{Corollary}
\newtheorem{proposition}[theorem]{Proposition}

{
\theoremstyle{definition}
\newtheorem{definition}{Definition}
\newtheorem{example}{Example}
\newtheorem{remark}{Remark}
}

\newcommand{\Dec}{\text{Dec}}
\newcommand{\floorsig}{\lfloor \log\frac{1}{\sigma} \rfloor}
\newcommand{\floorn}{\lfloor \log n \rfloor}
\newcommand{\tp}[1]{\theta^{(#1)}}
\newcommand{\htp}[1]{\hat\theta^{(#1)}}

\newcommand{\Xp}[1]{X^{(#1)}}

\newcommand{\bp}[1]{b^{(#1)}}

\newcommand{\papertitle}{Distributed Gaussian Mean Estimation under Communication Constraints: Optimal Rates and Communication-Efficient Algorithms}

\arxiv{0000.0000} 

\begin{document} 
	
\begin{frontmatter}
	
\title{ \papertitle \thanksref{T1}}  
\runtitle{Distributed Gaussian Mean Estimation}
\thankstext{T1}{The research was supported in part by NSF Grant DMS-1712735 and NIH grants R01-GM129781 and R01-GM123056.}
\begin{aug}
	\author{\fnms{T. Tony} \snm{Cai}\ead[label=e1]{tcai@wharton.upenn.edu}}
	\and
	\author{\fnms{Hongji} \snm{Wei}\ead[label=e2]{hongjiw@wharton.upenn.edu}
		\ead[label=u1,url]{http://www-stat.wharton.upenn.edu/$\sim$tcai/}}
	\runauthor{T. T. Cai and H. Wei}
	\affiliation{University of Pennsylvania}
	\address{DEPARTMENT OF STATISTICS\\
		THE WHARTON SCHOOL\\
		UNIVERSITY OF PENNSYLVANIA\\
		PHILADELPHIA, PENNSYLVANIA 19104\\
		USA\\
		\printead{e1}\\
		\phantom{E-mail:\ }\printead*{e2}\\
		\printead{u1}\\
	}
	
\end{aug}

\begin{abstract}
We study distributed estimation of a  Gaussian mean under communication constraints in a decision theoretical framework. 
Minimax rates of convergence, which characterize the tradeoff between the communication costs and statistical accuracy,  are established in both the univariate and multivariate settings. Communication-efficient and statistically optimal procedures are developed. In the univariate case, the optimal rate depends only on the total communication budget, so long as each local machine has at least one bit. However, in the multivariate case,  the minimax rate depends on the specific allocations of the communication budgets among the local machines. 

Although optimal estimation of a Gaussian mean is relatively simple in the conventional setting, it is quite involved under the communication constraints, both in terms of the optimal procedure design and lower bound argument.  The techniques developed in this paper can be of independent interest.  An essential step is the decomposition of the minimax estimation problem into two stages, localization and refinement. This critical decomposition provides a framework for both the lower bound analysis and optimal procedure design.

\end{abstract}



\end{frontmatter}


\section{Introduction}
    
In the conventional statistical decision theoretical framework, the focus is on the centralized setting where all the data are collected together and directly available. The main goal  is to develop  optimal (estimation,  testing, detection, ...) procedures, where optimality is understood with respect to the sample size and parameter space. Communication/computational costs are not part of the consideration. 

In the age of big data, communication/computational concerns associated with a statistical procedure are becoming increasingly important in contemporary applications. One of the difficulties for analyzing large datasets is that data are distributed, instead of in a single centralized location. This setting arises naturally in many statistical practices.
    
    \begin{itemize}
        \item {\bf Large datasets.} When the datasets are too large to be stored on a single computer or data center, it is natural to divide the whole dataset into multiple computers or data centers, each assigned a smaller subset of the full dataset. Such is the case for a wide range of applications.
        \item {\bf Privacy and security.} Privacy and security concerns can also cause the decentralization of the datasets. For example, medical and financial institutions often collect datasets that contain sensitive and valuable information. For privacy and security  reasons, the data cannot be released to a third party for a centralized analysis and need to be stored in different and secure places while performing data analysis. 
    \end{itemize}
    

Learning from distributed datasets, which is called {\it distributed learning}, has attracted much recent attention. For example, Google AI proposed a machine learning setting called ``Federated Learning" \citep{mcmahan2017federated}, which develops a high-quality centralized model while the training data remain distributed over a large number of clients.  Figure \ref{pic_dist_computers} provides a simple illustration of a distributed learning network. In addition to advances on architecture design for distributed learning in practice, there is also an increasing amount of literature on distributed learning theories, including \cite{jordan2018communication}, \cite{battey2018distributed},  \cite{dobriban2018distributed}, and \cite{fan2019distributed} in statistics, computer science, and information theory communities. Several distributed learning procedures with some theoretical properties have been developed in recent works. However, they do not impose any communication constraints on the proposed procedures thus fail to characterize the relationship between the communication costs and statistical accuracy. Indeed, in a decision theoretical framework,  if no communication constraints are imposed, one can always output the original data from the local machines to the central machine and treat the problem same as in the conventional centralized setting.

The study on how the communication constraints compromise the estimation accuracy in the distributed settings has a long history. Dating back to 1980's, \cite{zhang1988estimation} proposed an asymptotic unbiased distributed estimator and calculated its variance. In recent years, there is an emerging literature focusing on distributed Gaussian mean estimation under the communication constraints. \cite{garg2014communication} provided a bound on the bits of communication needed to achieve the centralized minimax risk. \cite{zhang2013information, braverman2016communication} introduced information-theoretical tools to prove lower bounds on the minimax rate for Gaussian mean estimation under communication constraints. \cite{han2018geometric} developed a geometric lower bound for distributed Gaussian mean estimation. Other similar settings and distribution families were also studied in \cite{luo2005universal,pmlr-v80-zhu18a, kipnis2019mean, hadar2019distributed, szabo2019asymptotic}.   
        
For large-scale data analysis,  communications between machines can be slow and expensive and limitation on bandwidth and communication sometimes becomes the main bottleneck on statistical efficiency.  It is therefore necessary to take communication constraints into consideration when constructing statistical procedures. When the communication  budget is limited, the algorithm must carefully ``compress" the information contained in the data as efficiently as possible, leading to a trade-off between communication costs and statistical accuracy.  The precisely quantification of this trade-off is an important and challenging problem.
    
Estimation of a Gaussian mean occupies a central position in parametric statistical inference.  In the present paper we consider  distributed Gaussian mean estimation under the communication constraints  in both the univariate and multivariate settings. Although optimal estimation of a Gaussian mean is a relatively simple problem in the conventional setting, this problem is quite involved under the communication constraints, both in terms of the construction of the rate optimal distributed estimator and the lower bound argument. Optimal distributed estimation of a Gaussian mean also serves as a starting point for investigating other more complicated statistical problems in distributed learning including distributed high-dimensional linear regression and distributed large-scale multiple testing.

    \subsection{Problem formulation}
        
 We begin by giving a formal definition of {\bf transcript},  \textbf{distributed estimator}, and \textbf{distributed protocol}.   
        Let $\mathcal{P} = \{P_{\theta} : \theta \in \Theta \}$ be a parametric family of distributions supported on space $\mathcal{X}$, where $\theta\in \Theta\subseteq  {\mathbb R}^d$ is the parameter of interest. Suppose there are $m$ local machines and a central machine, where the the local machines contain the observations and the central machine produces the final estimator of $\theta$ under the communication constraints between the local and central machines. More precisely, suppose we observe  i.i.d. random samples drawn from a distribution $P_{\theta}\in \mathcal{P}$:
\[  X_i \iid P_{\theta}, \quad i=1, \ldots, m, \]
where the $i$-th local machine has access to $X_i$ only. 

For $i=1, ..., m$, let  $b_i\ge 1$ be a positive integer and the  $i$-th local machine can only transmit $b_i$ bits to the central machine. That is, the observation $X_i$ on the $i$-th local machine needs to be processed to a binary string of length $b_i$ by a (possibly random) function $\Pi_i: \mathcal{X} \to \{0,1\}^{b_i}$. The resulting string $Z_i \triangleq \Pi_i(X_i)$, which is called the \textbf{transcript} from the $i$-th machine, is then transmitted to the central machine. Finally, a \textbf{distributed estimator} $\hat \theta$ is constructed on the central machine based on the transcripts $Z_1, Z_2, ..., Z_m$, 
        \[ \hat\theta = \hat\theta(Z_1, Z_2, ..., Z_m). \]
        The above scheme to obtain a distributed estimator $\hat\theta$  is called a \textbf{distributed protocol}. The class of distributed protocols with communication budgets $b_1,b_2,...,b_m$ is defined as
        \begin{align*}
            \mathcal{A}(b_1,b_2,...,b_m) = \{(\hat\theta, \Pi_1, \Pi_2,...,\Pi_m): \; &\Pi_i: \mathcal{X} \to \{0,1\}^{b_i}, \; i=1,2,...,m, \\
            &\hat\theta = \hat\theta(\Pi_1(X_1), ..., \Pi_m(X_m))\}. 
        \end{align*}   
        
               \begin{figure}[htb]
       		\begin{subfigure}[t]{0.54\textwidth}
       			\includegraphics[width=2.6in, height=1.9in]{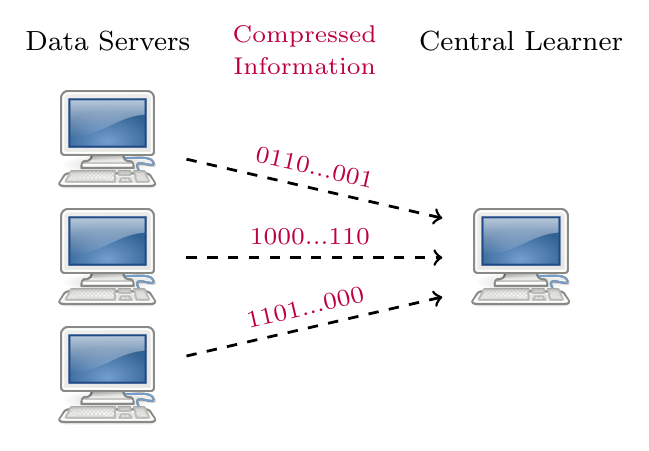}
       			\caption{\small Distributed learning network}
       			\label{pic_dist_computers}
       		\end{subfigure}
        	\begin{subfigure}[t]{0.44\textwidth}
        		\includegraphics[width=0.9\textwidth]{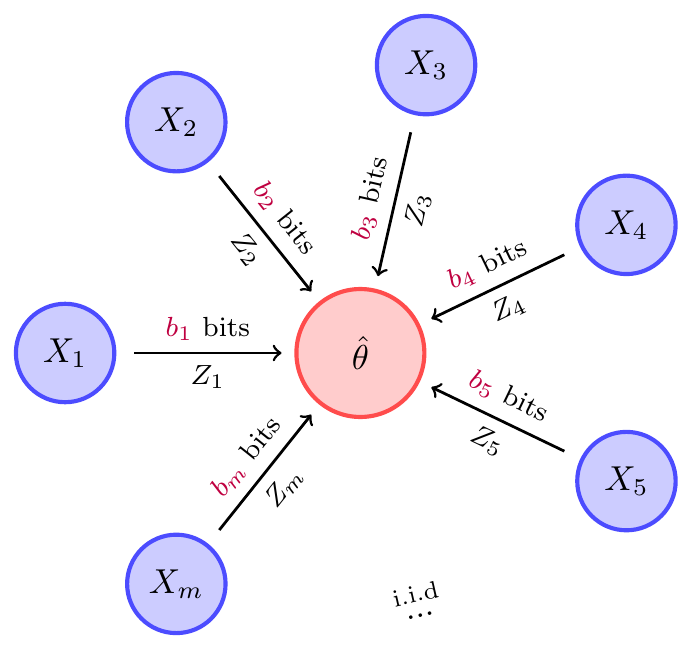}
        		\caption{\small Distributed protocol}
        		\label{pic_dist_machines}
        	\end{subfigure}
        	\caption{\small (a) Left panel: An illustration of a distributed learning network. Communication between the data servers and the central learner is necessary in order to learn from distributed datasets.  (b) Right panel: An illustration of distributed protocol. The $i$-th machine can only transmit a $b_i$ bits transcript to the central machine. } \vspace{-8pt}
        \end{figure}
        
        We use $b_{1:m}$ as a shorthand for $(b_1,b_2,...,b_m)$ and denote $\hat\theta \in \mathcal{A}(b_{1:m})$ for $(\hat\theta, \Pi_1, ..., \Pi_m) \in \mathcal{A}(b_{1:m})$. We shall always assume $b_i \geq 1$ for all $i=1,2,...,m$, i.e. each local machine can transmit at least one bit to the central machine. Otherwise, if no communication is allowed from any of the local machines, one can just exclude  those local machines and treat the problem as if there are fewer local machines available. Figure \ref{pic_dist_machines} gives  a simple illustration for the distributed protocols.
        
As usual, the estimation accuracy of a distributed estimator $\hat\theta$ is measured by the mean squared error (MSE), $\E_{P_{\theta}} \Vert \hat\theta - \theta \Vert_2^2$, 
where the expectation is taken over the randomness in both the data and construction of the transcripts and estimator. 
 As in the conventional decision theoretical framework, a quantity of particular interest in distributed learning is the minimax risk for the distributed protocols 
        \[  \inf_{\hat\theta\in \mathcal{A}(b_{1:m})} \sup_{P_{\theta} \in \mathcal{P}} \E_{P_{\theta}} \Vert \hat\theta - \theta \Vert_2^2,    \] 
which characterizes the difficulty of the distributed learning problem under the communication constraints $b_{1:m}$.
As mentioned earlier,  in a rigorous decision theoretical formulation of distributed learning, the communication constraints are essential. Without the constraints, one can always output the original data from the local machines to the central machine and the problem is then reduced to the usual centralized setting.

    \subsection{Distributed estimation of a univariate Gaussian mean}
            
We first consider distributed estimation of a univariate Gaussian mean under the communication constraints $b_{1:m}$, where $P_{\theta} =  N(\theta,\sigma^2)$ with $\theta\in [0, 1]$ and the variance $\sigma^2$ known. Note that by a sufficiency argument, the case where each local machine has access to $k$ i.i.d. samples from $N(\theta,\sigma^2)$ is the same. 
        
Our analysis in Section \ref{sec.univ} establishes the following minimax rate of convergence for distributed univariate Gaussian mean estimation under the communication constraints $b_{1:m}$,
\beq
\inf_{\hat\theta\in\mathcal{A}(b_{1:m})}\sup_{\theta\in [0,1]} \E(\hat\theta - \theta)^2 \asymp \begin{cases} 2^{-2B} &\text{ if } B< \log \frac{1}{\sigma} + 2 \\ \frac{\sigma^2}{(B-\log \frac{1}{\sigma})} &\text{ if } \log \frac{1}{\sigma} + 2 \leq B < \log \frac{1}{\sigma} + m \\ \min\left\{\frac{\sigma^2}{m}, 1 \right\} &\text{ if } B \geq \log \frac{1}{\sigma} + m \end{cases},
\label{intro_univ_rate} 
\eeq       
where $B = \sum_{i=1}^m b_i$ is the total communication budgets, and $a \asymp b$ denotes $cb \leq a \leq Cb$ for some constants $c, C>0$.

\begin{figure}[htb]
   	\centering
   	\includegraphics[width=3.5in]{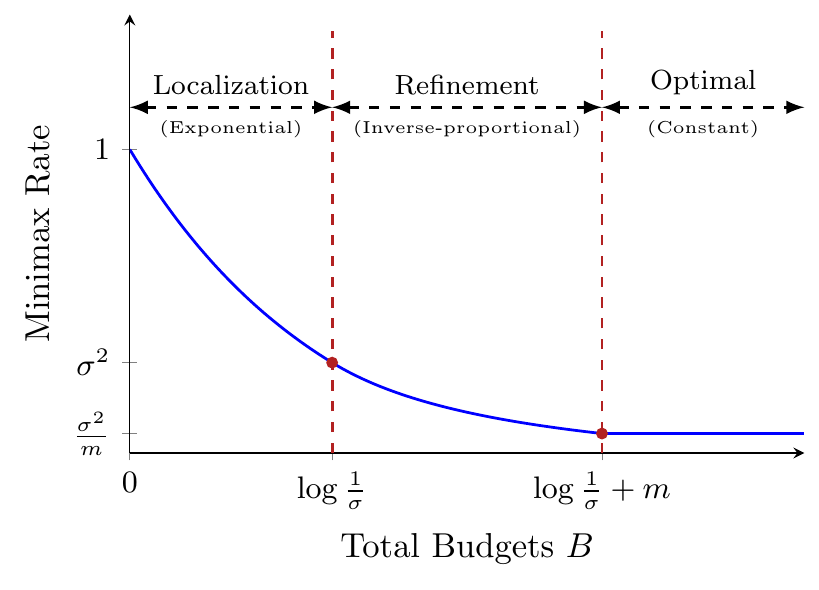}
   	\caption{\small The minimax rate of univariate Gaussian mean estimation under communication constraints has 3 phases: localization, refinement and optimal-rate.} \vspace{-8pt}
   	\label{pic_univ_rate} 	
   \end{figure}

   The above minimax rate characterizes the trade-off between the communication costs and statistical accuracy for univariate Gaussian mean estimation. An illustration of the minimax rate is shown in Figure \ref{pic_univ_rate}.
   
        The minimax rate \eqref{intro_univ_rate} is interesting in several aspects. First, the optimal rate of convergence only depends on the total communication budgets $B = \sum_{i=1}^m b_i$, but not the specific allocation of the communication budgets among the $m$ local machines, as long as each machine has at least one bit. 
        Second, the rate of convergence has three different phases:          
        \begin{enumerate}
            \item Localization phase. When $B< \log_2\frac{1}{\sigma}+2$, as a function of $B$, the minimax risk decreases  fast at an exponential rate. In this phase, having more communication budget is very beneficial in terms of improving the estimation accuracy.    
                    
            \item Refinement phase. When $\log_2\frac{1}{\sigma}+2 \leq B < \log_2\frac{1}{\sigma}+m$, as a function of $B$, the minimax risk decreases relatively slowly and is inverse-proportional to the total communication budget $B$.
            
            \item Optimal-rate phase. When $B \geq \log_2\frac{1}{\sigma}+m$, the minimax rate does not depend on $B$, and is the same as in the centralized setting where all the data are combined \citep{bickel1981minimax}.
        \end{enumerate}
        
         An essential technique for solving this problem is the decomposition of the minimax estimation problem into two steps, {\it localization} and {\it refinement}. This critical decomposition provides a framework for both the lower bound analysis and optimal procedure design. In the lower bound analysis, the statistical error is decomposed into ``localization error" and ``refinement error". It is shown that one of these two terms is inevitably large under the communication constraints. In our optimal procedure called MODGAME, bits of the transcripts are divided into three types: crude localization bits, finer localization bits, and refinement bits. They compress the local data in a way that both the localization and refinement errors can be optimally reduced. Further technical details and discussion are presented in Section \ref{sec.univ}.

    \subsection{Distributed estimation of a Multivariate Gaussian mean}
            
We then consider the multivariate case under the communication constraints $b_{1:m}$, where $P_{\theta} =  N_d(\theta,\sigma^2I_d)$ with $\theta\in [0, 1]^d$ and the noise level $\sigma$ is known. Similar to the univariate case,  the goal is to optimally estimate the mean vector $\theta$ under the squared error loss. 
        
        The construction and the analysis given in Section \ref{sec.multi} show that the minimax rate of convergence in this case is given by
        \beq  \inf_{\hat\theta\in\mathcal{A}(b_{1:m})}\sup_{\theta\in [0,1]^d} \E \Vert \hat\theta - \theta \Vert_2^2 \asymp \begin{cases} 2^{-2B/d}d &\text{ if } B/d < \log_2 \frac{1}{\sigma} + 2 \\ \frac{d\sigma^2}{(B/d-\log_2 \frac{1}{\sigma})} &\text{ if } \log_2\frac{1}{\sigma} + 2 \leq B/d < \log_2 \frac{1}{\sigma} + \max\{m',2\} \\ d\min\left\{\frac{\sigma^2}{m'},1\right\} &\text{ if } B/d \geq \log_2 \frac{1}{\sigma} + \max\{m',2\} \end{cases} \label{intro_mutli_rate} \eeq
        where $B = \sum_{i=1}^m b_i$ is the total communication budgets and $m' = \sum_{i=1}^m \min\left\{\frac{b_i}{d}, 1\right\}$ is the ``effective sample size".
        
        The minimax rate in the multivariate case \eqref{intro_mutli_rate} is an extension of its univariate counterpart \eqref{intro_univ_rate}, but it also has its distinct features, both in terms of the estimation procedure and lower bound argument. Intuitively, the total communication budgets $B$ are evenly divided into $d$ parts so that roughly $B/d$ bits can be used to estimate each coordinate. Because there are $d$ coordinates, the risk is multiplied by $d$. The effective sample size $m'$ is a special and interesting quantity in multivariate Gaussian mean estimation. This quantity suggests that even when the total communication budgets are sufficient, the rate of convergence must be larger than the benchmark $d\min\left\{\frac{\sigma^2}{m'},1\right\}$. There is a gap between the distributed optimal rate and centralized optimal rate if $m' \ll m$. See Section \ref{sec.multi} for further technical details and discussion.
                
Although the interplay between communication costs and statistical accuracy has drawn increasing recent attention, to the best of our knowledge, the present paper is the first to establish a sharp minimax rate for distributed Gaussian mean estimation. Compared to our results, none of the previous results turns out to be sharp in general.  The techniques developed in this paper, both for the lower bound analysis and construction of the rate optimal procedure, can be of independent interest. Our lower bound argument was inspired by the earlier work on the strong data processing inequality proposed in \cite{zhang2013information, braverman2016communication, raginsky2016strong}. 

    \subsection{Organization of the paper}
    
We finish this section with notation and definitions that will be used in the rest of the paper. Section \ref{sec.univ} studies distributed estimation of a univariate Gaussian mean under communication constraints and Section \ref{sec.multi} considers the multivariate case.  The numerical performance of the proposed distributed estimators is investigated in Section \ref{sec.simulation}  and further research directions are discussed in Section \ref{sec.discussion}. For reasons of space, we prove the main results for the univariate case in Section  \ref{sec.proof} and defer the proofs of the results for the multivariate case and technical lemmas to the Supplementary Material \citep{CaiWei2019Supplement}.

    \subsection{Notation and definitions}
    For any $a \in \reals$, let $\lfloor a \rfloor$ denote the floor function (the largest integer not larger than $a$). Unless otherwise stated, we shorthand $\log a$ as the base 2 logarithmic of $a$. For any $a,b \in \reals$, let $a\wedge b \triangleq \min\{a,b\}$ and $a\vee b \triangleq \max\{a,b\}$. For any vector $a$, we will use $a^{(k)}$ to denote the $k$-th coordinate of $a$, and denote by $\Vert a \Vert \triangleq \sqrt{\sum_{k} \left(a^{(k)}\right)^2}$  its $l_2$ norm. For any set $S$, let $S^k \triangleq S \times S \times ... \times S$ be the Cartesian product of $k$ copies of $S$. Let $\ind{\cdot}$ denote the indicator function taking values in $\{0,1\}$. 
    
    For any discrete random variables $X, Y$ supported on $\mathcal{X}, \mathcal{Y}$, the entropy $H(X)$, conditional entropy $H(X|Y)$, and mutual information $I(X;Y)$ are defined as
    \beas 
    H(X) &\triangleq& - \sum_{x \in \mathcal{X}} \P(X=x)\log\P(X=x),    \\
    H(X|Y) &\triangleq& - \sum_{x \in \mathcal{X}, y \in \mathcal{Y}} \P(X=x, Y=y)\log\P(X=x|Y=y),   \\
    I(X; Y) &\triangleq& \sum_{x \in \mathcal{X}, y \in \mathcal{Y}} \P(X=x, Y=y)\log\frac{\P(X=x|Y=y)}{\P(X=x)}.   
    \eeas

  
\section{Distributed Univariate Gaussian Mean Estimation} 
\label{sec.univ} 

In this section we consider distributed estimation of a univariate Gaussian mean, where one observes  on $m$ local machines i.i.d. random samples:
\[  X_i \iid N(\theta,\sigma^2), \quad i=1, \ldots, m, \]
under the constraints that the $i$-th machine has access to $X_i$ only and can transmit $b_i$ bits only to the central machine. 
We denote by $\mathcal{P}_\sigma^1$ the  Gaussian location family
\[  \mathcal{P}_\sigma^1 = \left\{ N(\theta,\sigma^2) : \theta \in [0,1]  \right\},  \]
where $\theta \in [0,1]$ is the mean parameter of interest and the variance $\sigma^2$ is known.
For given communication budgets $b_{1:m}$ with $b_i \ge 1$ for $i=1, \ldots, m$, the goal is to optimally estimate the mean $\theta$ under the squared error loss. A particularly interesting quantity is the minimax risk under the communication constraints, i.e., the  minimax risk for the distributed protocol $\mathcal{A}(b_{1:m})$:
\[  R_1(b_{1:m}) = \inf_{\hat\theta\in\mathcal{A}(b_{1:m})}\sup_{\theta\in [0,1]} \E(\hat\theta - \theta)^2,\]
which characterizes the difficulty of the estimation problem under the communication constraints. We are also interested in constructing a computationally efficient algorithm that achieves the minimax optimal rate. 

We first introduce an estimation procedure and provide an upper bound for its performance and then establish a matching lower bound on the minimax rate. The upper and lower bounds together establish  the minimax rate of convergence and the optimality of the proposed estimator. 
 
\subsection{Estimation procedure - MODGAME}
\label{MODGAME.sec}

We begin with the construction of an estimation procedure under the communication constraints and provide a theoretical analysis of the proposed procedure. The procedure, called MODGAME (Minimax Optimal Distributed GAussian Mean Estimation), is a deterministic procedure that generates a distributed estimator $\hat\theta_D$ under the distributed protocol $\mathcal{A}(b_{1:m})$. We divide the discussion into two cases: $\sigma < 1$ and $\sigma \geq 1$. 

\subsubsection{MODGAME procedure when $\sigma < 1$}
\label{sigma<1.sec}

When $\sigma < 1$, MODGAME consists of two steps: localization and refinement. Roughly speaking, the first step utilizes $\log\frac{1}{\sigma} + o(B-\log\frac{1}{\sigma})$ bits, out of  the total budget $B=\sum_{i=1}^m b_i$ bits, for localization to roughly locate where $\theta$ is, up to $O(\sigma)$ error. Building on the location information,  the remaining $B-\log\frac{1}{\sigma}$ bits are used for refinement to further increase the accuracy of the estimator. Detailed theoretical analysis will show that the optimality of the final estimator.

Before describing the MODGAME procedure in detail, we define several useful functions that will be used to generate the transcripts. For any interval $[L,R]$, let $\tau_{[L,R]}: \reals \to [L,R]$ be the truncation function defined by
\beq \tau_{[L,R]}(x) = \begin{cases} L \quad {\rm if} \quad x\leq L  \\ x \quad {\rm if} \quad L<x<R \\ R \quad {\rm if} \quad x\geq R \end{cases}.  \label{def_tau} \eeq

For any integer $k\geq 0$, denote $g_k: \reals \to \{0,1\}$ be the $k$-th Gray function defined by
\[ g_k(x) \triangleq \begin{cases} 0 \quad {\rm if} \quad \lfloor 2^k\tau_{[0,1]}(x) \rfloor \  { \rm mod } \  4 = 0 \  {\rm or }\  3 \\ 1 \quad {\rm if} \quad \lfloor 2^k\tau_{[0,1]}(x) \rfloor \  { \rm mod } \  4 = 1 \  {\rm or }\  2 .  \end{cases}   \]

Similarly we denote by $\bar{g}_k: \reals \to \{0,1\}$ the $k$-th conjugate Gray function defined by
\[ \bar{g}_k(x) \triangleq \begin{cases} 0 \quad {\rm if} \quad \lfloor 2^k\tau_{[0,1]}(x) \rfloor \  { \rm mod } \  4 = 0 \  {\rm or }\  1 \\ 1 \quad {\rm if} \quad \lfloor 2^k\tau_{[0,1]}(x) \rfloor \  { \rm mod } \  4 = 2 \  {\rm or }\  3 .  \end{cases}   \]

To unify the notation we set $g_k(x) \equiv \bar{g}_k(x) \equiv 0$ if $k < 0$.

\begin{figure}[htb]
	\centering
	\includegraphics[width=0.6\textwidth]{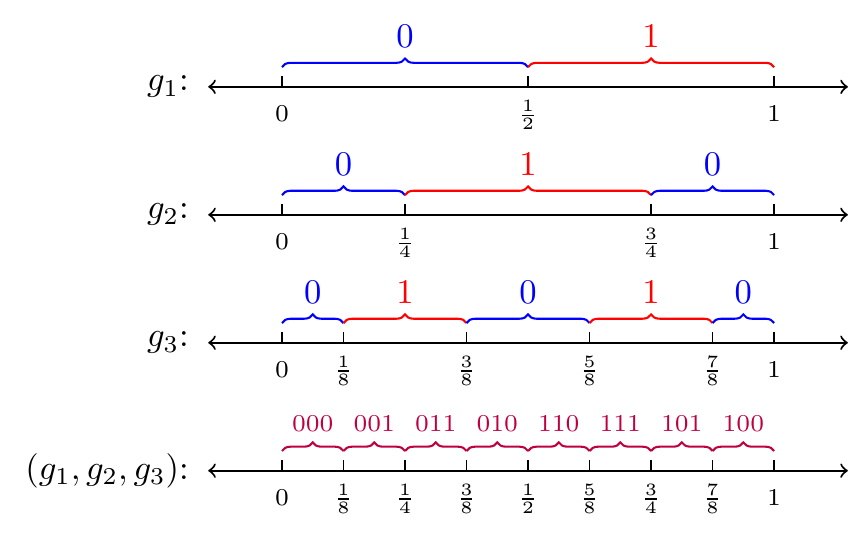}
	\caption{\small An illustration of the Gray functions and Gray codes.}
	\vspace{-15pt}
	\label{pic_gray_ex}
\end{figure}

It is worth mentioning that these Gray functions mimic the behavior of the Gray codes (for reference see \cite{savage1997survey}). Fix $K\geq 1$, if we treat $(g_1(x), g_2(x), ..., g_K(x))$ as a string of code for any source $x\in[0,1]$, then those $x$ within the interval $[2^{-K}(s-1), 2^{-K}s)$ where $s$ is a integer will match the same code. Moreover, the code for adjacent intervals only differs by one bit, which is also a key feature for the Gray codes. This key feature guarantees the robustness of the Gray codes. Such robustness makes the Gray functions very useful for distributed estimation. An example for $K=3$ is shown in Figure \ref{pic_gray_ex} to better illustrate the behavior of the Gray functions.

Define the refinement function $h(x): \reals \to \{0,1\}$ and the conjugate refinement function $\bar{h}(x): \reals \to \{0,1\}$ by
\[  h(x) \triangleq  \lfloor 2^{(\floorsig-7)} x \rfloor \  { \rm mod } \  2  \quad {\rm and} \quad  \bar{h}(x) \triangleq \lfloor 2^{(\floorsig-7)} x - \frac{1}{2} \rfloor \  { \rm mod } \  2 . \]

For any function $f$, define the convolution function
\[ \Phi_f(x) \triangleq \E_{X\sim N(x,\sigma^2)} f(X) = \int_{-\infty}^\infty \frac{1}{\sqrt{2\pi}\sigma}e^{-\frac{(y-x)^2}{2\sigma^2}}f(y)dy. \]


For any $K \geq 1$, let $\Dec_K(y_1,y_2,...,y_K): \{0,1\}^K \to 2^{[0,1]}$ be the decoding function defined by
\[ \Dec_K(y_1,y_2,...,y_n) \triangleq \{x\in[0,1] : g_k(x) = y_k \quad {\rm for} \ k=1,2,...,K    \} .  \]

Last, we define the distance between a point $x \in \reals$ and a set $S \subseteq \reals$ as
\[ d(x, S) \triangleq \min_{y \in S} |x-y| . \]

We are now ready to introduce the MODGAME procedure in detail. Again, we divide into three cases.

\medskip\noindent
\textbf{Case 1: $B < \log\frac{1}{\sigma} + 2$.}
In this case, the output is the values of the first $B$ localization bits from local machines, where the $k$-th localization bit is defined as the value of the function $g_k(\cdot)$ evaluated on the local sample. The procedure can be described as follows.

\begin{enumerate}
    \item[\textbf{Step 1:}] {\it Generate transcripts on local machines.}
    Define $s_0 = 0$ and $s_i = \sum_{j=1}^{i} b_j$ for  $i=1, \ldots, m$. On the $i$-th machine, the transcript $Z_i$ is concatenated by the $(s_{i-1} + 1)$-th, $(s_{i-1} + 2)$-th, ..., $(s_{i-1} + b_i)$-th Gray functions evaluated at $X_i$. That is,
    \[ Z_i = (U_{s_{i-1} + 1}, U_{s_{i-1} + 2}, ..., U_{s_{i-1} + b_i}) ,  \]
    where $U_{s_{i-1} + k} \triangleq g_{s_{i-1} + k}(X_i)$ for $k = 1,2,...,b_i$. 
    
    \item[\textbf{Step 2:}] {\it Construct distributed estimator $\hat\theta_D$.}    
    Now we collect the bits $U_1,U_2,...,U_B$ from the transcripts $Z_1,Z_2,...,Z_m$. Note that $Z_k$ is the $k$-th Gray function evaluate at a random sample drawn from $N(\theta,\sigma^2)$, one may reasonably "guess" that $U_k \approx g_k(\theta)$. By this intuition, we set $\hat\theta_D$ to be the minimum number in the interval $\Dec_B(U_1,U_2,...,U_B)$, i.e.
    \[ \hat\theta_D = \min\{x: x\in\Dec_B(U_1,U_2,...,U_B) \} .  \]
\end{enumerate}

\medskip\noindent
\textbf{Case 2:} $\log\frac{1}{\sigma} + 2 \leq B \leq \log\frac{1}{\sigma} + m$.
Let
\beq   n \triangleq \max\left\{ s \in \naturals : \lfloor \log s \rfloor^2 + 2s \leq B - \lfloor\log\frac{1}{\sigma}\rfloor \right\},  \label{def_n}  \eeq
and define finer localization functions:
\begin{equation}
    \begin{aligned}
        f_1(x) &\triangleq g_{\floorsig - \floorn - 2 }(x) ,   \\ 
        f_2(x) &\triangleq \bar{g}_{\floorsig - \floorn - 2 }(x) ,    \\
        f_k(x) &\triangleq g_{\floorsig - \floorn - 4 + k }(x) \text{ for } k\geq 3 . \\
    \end{aligned}
    \label{def_f}
\end{equation}  

In this case the total communication budget is divided into 3 parts: crude localization bits (roughly $\lfloor\log\frac{1}{\sigma}\rfloor$ bits), finer localization bits ($\lfloor \log n \rfloor^2$ bits), and refinement bits ($2n$ bits). 
The crude localization bits are the values of the functions $g_1(\cdot), g_2(\cdot), ..., g_{\lfloor\log\frac{1}{\sigma}\rfloor}(\cdot)$, each evaluated on a local sample. We denote those resulting binary bits by $U_1, U_2, ..., U_{\lfloor\log\frac{1}{\sigma}\rfloor}$.
The finer localization bits are the values of the functions $f_1(\cdot), f_2(\cdot), ..., f_{\floorn}(\cdot)$, each function is evaluated on $\lfloor \log n \rfloor$ different local samples. The function values of $f_k(\cdot)$ are denoted by $W_{k,1}, W_{k,2}, ..., W_{k,n}$.
The refinement bits are the values of the function $h(\cdot)$, evaluated on $n$ local samples; and the values of the function $\bar{h}(\cdot)$, evaluated on $n$ different local samples. The resulting binary bits are denoted by $V_1,V_2,...,V_{n}$ and $\bar{V}_1,\bar{V}_2,...,\bar{V}_{n}$ respectively.

These three types of bits are assigned to local machines by the following way: (1) Among all $m$ machines, there are $\floorn^2$ local machines who will output transcript consist of 1 finer localization bit and $b_i-1$ crude localization bits. (2) Among all $m$ machines, there are $2n$ local machines who will output transcript consist of 1 refinement bit and $b_i-1$ crude localization bits. (3) The remain $m-(\floorn^2 + 2n)$ machines will output transcript consist of $b_i$ crude localization bits. The above assignment is feasible because
\[  \floorn^2 + 2n \leq B - \floorsig \leq m.  \]

It is worth mentioning that every finer localization bits and every refinement bits are assigned to different machines. Intuitively, this is because we need these bits to be independent so that we can gain enough information for estimation. 
See Figure \ref{pic_modgame} for an overview of the MODGAME procedure. The procedure can be summarized as follows:

\begin{figure}[t]
	\centering
	\includegraphics[width=0.96\textwidth]{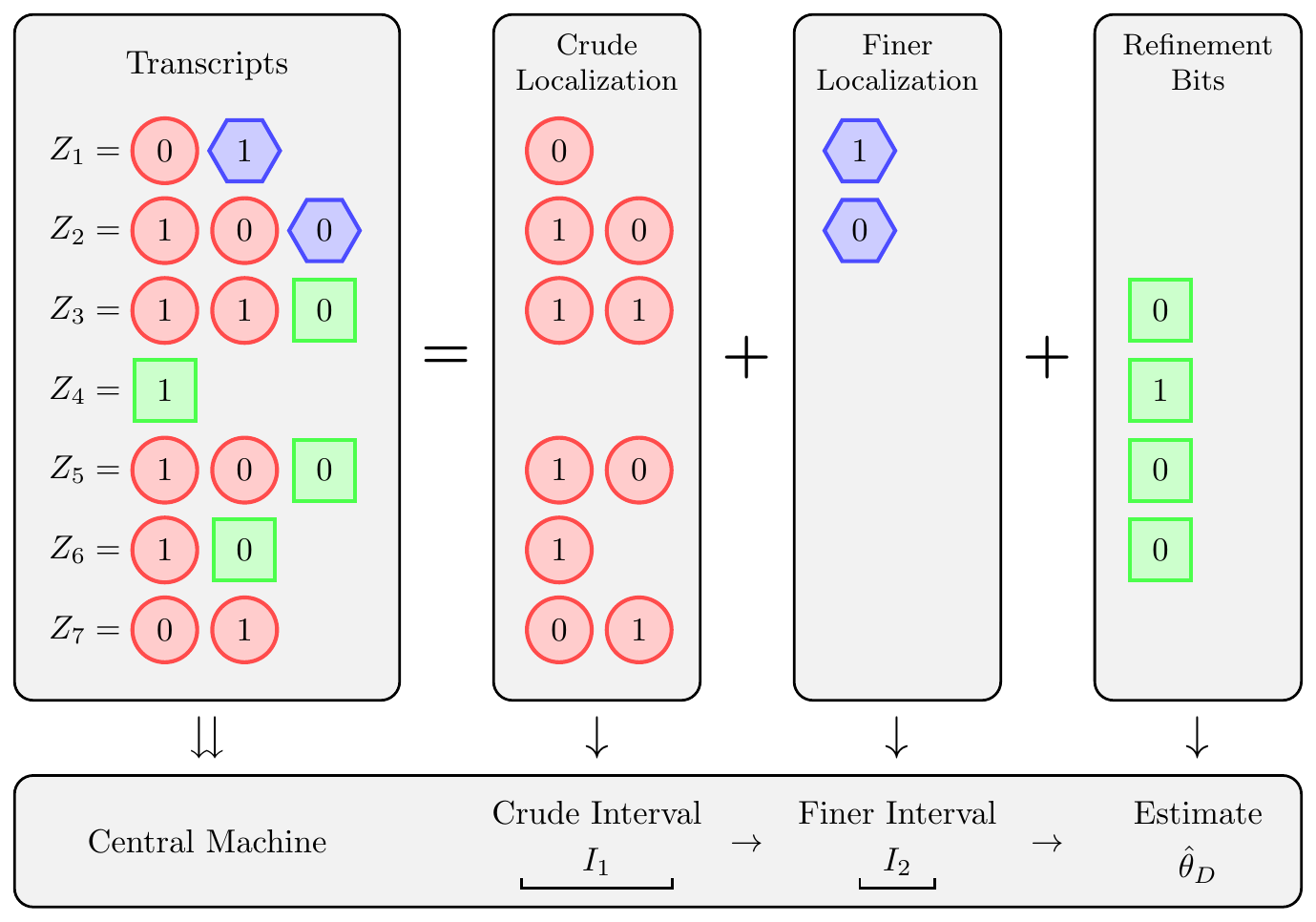}
	\caption{\small An illustration of MODGAME. The bits in the transcripts are transmitted to the central machine and divided into three types: crude localization bits, finer localization bits, and refinement bits. One then constructs on the central machine a crude interval $I_1$, a finer interval $I_2$, and the final estimate $\hat\theta_D$ step by step.}
	\vspace{-15pt}
	\label{pic_modgame}
\end{figure}

\begin{enumerate}
    \item[\textbf{Step 1:}] {\it Generate transcripts on local machines.}
    Define $s_i = \sum_{j=1}^{i} (b_j - \ind{j\leq \lfloor \log n \rfloor^2 + 2n})$ and $s_0 = 0$. On the $i$-th machine:
    \begin{itemize}
    	\item If $(j-1)\floorn + 1 \leq i \leq j \floorn$ for some integer $1\leq j \leq \floorn$, output 
    	\[  Z_i = (U_{s_{i-1}+1}, U_{s_{i-1}+2},...,U_{s_{i-1}+b_i-1}, W_{j, i-(j-1)\floorn}) ;  \]
    	(If $b_i=1$, just output $Z_i = W_{j, i-(j-1)\floorn}$.) 
    	
    	\item If $\floorn^2 + 1 \leq i \leq \floorn^2 + n$, output
    	\[  Z_i = (U_{s_{i-1}+1}, U_{s_{i-1}+2},...,U_{s_{i-1}+b_i-1},V_{i-\floorn^2})  ;  \]
    	(If $b_i=1$, just output $Z_i = V_{i-\floorn^2}$.) 
    	
    	\item If $\floorn^2 + n + 1 \leq i \leq \floorn^2 + 2n$, output
    	\[  Z_i = (U_{s_{i-1}+1}, U_{s_{i-1}+2},...,U_{s_{i-1}+b_i-1},\bar{V}_{i-(\floorn^2 + n)})  ;  \]
    	(If $b_i=1$, just output $Z_i = \bar{V}_{i-(\floorn^2 + n)}$.) 
    	
    	\item If $i \geq \floorn^2 + 2n + 1$, output
    	\[  Z_i = (U_{s_{i-1}+1}, U_{s_{i-1}+2},...,U_{s_{i-1}+b_i}) .   \]
    \end{itemize}

	where the above binary bits are calculated by
	\beas U_{s_{i-1}+k} \triangleq g_{s_{i-1}+k}(X_i)  &{\rm for} &i = 1,2,...,m; \  k=1,2,...,b_i. \\
	W_{j, i-(j-1)\floorn} \triangleq f_j(X_i)  &{\rm for} &j=1,2,...,\floorn - 1; \\ && i = (j-1)\floorn + 1, (j-1)\floorn + 2,...,j\floorn. \\
	V_{i-\floorn^2} \triangleq h(X_i)  &{\rm for}  &i = \floorn^2 + 1, \floorn^2 + 2, ..., \floorn^2 + n. \\
	 \bar{V}_{i-(\floorn^2 + n)} \triangleq \bar{h}(X_i)   &{\rm for}  &i = \floorn^2 + n +1, ..., \floorn^2 + 2n . \eeas
	
	\item[\textbf{Step 2:}] {\it Construct distributed estimator $\hat\theta_D$.}
	From transcripts $Z_1,Z_2,...,Z_m$, we can collect (a) crude localization bits $U_1,U_2,...,U_{\floorsig}$; (b) finer localization bits $W_{1,1}, W_{1,2}, ..., W_{\floorn, \floorn}$; (c) refinement bits $V_1,V_2,...,V_n$ and $\bar{V}_1, \bar{V}_2, ..., \bar{V}_n$. 
	
	\begin{enumerate}
		\item[\textbf{Step 2.1:}] First, we use crude localization bits $U_1, U_2, ,,, U_{\floorsig - \floorn - 3}$ to roughly locate $\theta$. The ``crude interval" $I_1$ will be obtained in this step. 
		
		(a) If $\floorsig - \floorn \leq 3$, just set $I_1 = I_1' = [0,1]$.
		
		(b) If $\floorsig - \floorn \geq 4$, let
		\beq  I_1' \triangleq \Dec_{\floorsig-\floorn-3}(U_1,U_2,...,U_{\floorsig-\floorn-3}) . \label{def_I1p} \eeq
	
		Then we further stretch $I_1'$ to a larger interval $I_1$ so that $I_1$ will double the length of $I_1'$:
		\beq I_1 \triangleq \left\{ x: d(x, I_1') \leq 2^{-(\floorsig - \floorn -2)} \right\} .
		\label{def_I1} \eeq
		
		\item[\textbf{Step 2.2:}] Then, we use finer localization bits to locate $\theta$ to a smaller interval of length roughly $O(\sigma)$. A "finer interval" $I_2$ will be generated in this step. For any $1\leq k \leq \floorn$, let
		\[  W_k = \ind{\sum_{j=1}^{\floorn} W_{k,j} \geq \frac{1}{2} \floorn }  \]
		be the majority voting summary statistic for $W_{k,1}, W_{k,2}, ..., W_{k,\lfloor \log n \rfloor}$. 
		
        (a) If $\floorsig - \floorn \leq 3$, and $\floorsig \leq 4$, let
        \[  I_2 = I_2' = [0,1] .  \]
        
    	(b) If $\floorsig - \floorn \leq 3$, and $\floorsig \geq 5$, let 
    	\beq I_2' \triangleq \Dec_{\floorsig - 4}(W_{\floorn-\floorsig+5},W_{\floorn-\floorsig+6}, ..., W_{\floorn}) . \label{def_I2p_trivialcase} \eeq
    	
    	Then we further stretch $I_2'$ to a larger interval $I_2$ so that $I_2$ will double the length of $I_2'$:
    	\[  I_2 \triangleq \left\{ x: d(x, I_2') \leq 2^{-(\floorsig - 3)} \right\} . \]
    	
    	(c) If $\floorsig - \floorn \geq 4$, let
    	\beq I_2' \triangleq  \left\{ x \in I_1: f_k(x) = W_k \text{ for all } k=1,2,...,\floorn \right\} \label{def_I2p} . \eeq
     
    	Lemma \ref{lm_ub_5} shows $I_2'$ is an interval. Then we further stretch $I_2'$ to a larger interval $I_2$ so that $I_2$ will double the length of $I_2'$:
    	\[  I_2 \triangleq \left\{ x: d(x, I_2') \leq 2^{-(\floorsig - 3)} \right\} . \]

    	\item[\textbf{Step 2.3:}] Finally, we use refinement bits $V_1, V_2, ..., V_n$ and $\bar{V}_1, \bar{V}_2, ..., \bar{V}_n$ to get an accurate estimate $\hat\theta_D$.  Lemma \ref{lm_ub_6} shows that one of the following two conditions must hold:
    	\[ I_2 \subseteq \left[(2j-\frac{3}{4})\cdot2^{-(\floorsig - 6)}, (2j+\frac{3}{4})\cdot2^{-(\floorsig - 6)}\right] \text{ for some } j\in \integers \] 
    	or 
    	\[ I_2 \subseteq \left[(2j+\frac{1}{4})\cdot2^{-(\floorsig - 6)}, (2j+\frac{7}{4})\cdot2^{-(\floorsig - 6)}\right] \text{ for some } j\in \integers . \]
    	So we can divide the procedure into the following two cases.
    	
    	(a) If $I_2 \subseteq [(2j-\frac{3}{4})\cdot2^{-(\floorsig - 6)}, (2j+\frac{3}{4})\cdot2^{-(\floorsig - 6)}] \text{ for some } j\in \integers$. Then $\Phi_h(x)$ is a strictly monotone function on $I_2$ (proved in Lemma \ref{lm_ub_6}). Denote 
    	\[  L_I \triangleq \inf_{x\in I_2}\Phi_h(x)  \quad {\rm and} \quad R_I \triangleq \sup_{x\in I_2}\Phi_h(x) . \] 
    	
    	By monotonicity, $\Phi_h$ is invertible on $I_2$. Let $\Phi_h^{-1}: [L_I, R_I] \to I_2$ be the inverse of $\Phi_h$, the distributed estimator $\hat\theta_D$ is given by
    	\beq \hat\theta_D = \Phi_h^{-1}\left(\tau_{[L_I,R_I]}\left(\frac{1}{n}\sum_{j=1}^n V_j\right)\right) \label{def_thetaD_1} \eeq
    	where $\tau_{[L_I,R_I]}$ is the truncation function defined in \eqref{def_tau}.
    	
    	(b) Otherwise, we have $I_2 \subseteq [(2j+\frac{1}{4})\cdot2^{-(\floorsig - 6)}, (2j+\frac{7}{4})\cdot2^{-(\floorsig - 6)}] \text{ for some } j\in \integers$. In this case $\Phi_{\bar{h}}(x)$ is a strictly monotone function on $I_2$ (proved in Lemma \ref{lm_ub_6}). Denote 
    	\[  \bar{L}_I \triangleq \inf_{x\in I_2}\Phi_{\bar{h}}(x)  \quad {\rm and} \quad \bar{R}_I \triangleq \sup_{x\in I_2}\Phi_{\bar{h}}(x) . \] 
    	
    	By monotonicity, $\Phi_{\bar{h}}$ is invertible on $I_2$. Let $\Phi_{\bar{h}}^{-1}: [\bar{L}_I, \bar{R}_I] \to I_2$ be the inverse of $\Phi_{\bar{h}}$, the distributed estimator $\hat\theta_D$ is given by
    	\beq \hat\theta_D = \Phi_h^{-1}\left(\tau_{[\bar{L}_I,\bar{R}_I]}\left(\frac{1}{n}\sum_{j=1}^n \bar{V}_j\right)\right) \label{def_thetaD_2} \eeq
    	where $\tau_{[\bar{L}_I,\bar{R}_I]}$ is the truncation function defined in \eqref{def_tau}.
	\end{enumerate}
    
\end{enumerate}

\medskip\noindent
\textbf{Case 3: $B > \log\frac{1}{\sigma} + m$.}
We only need to use part of the total communication budget $B$ as if we deal with the case $B = \lfloor \log\frac{1}{\sigma} \rfloor + m$. To be precise, we can always easily find $b_1', b_2', ..., b_m'$ so that $1\leq b_i' \leq b_i$ for  $i=1,2,...,m$ and 
\[  \sum_{i=1}^m b_i' = \lfloor \log\frac{1}{\sigma} \rfloor + m  . \]

Then we can implement the procedure introduced in Case 2 where we let the $i$-th machine only output a transcript of length $b_i'$.

\subsubsection{MODGAME procedure when $\sigma \geq 1$}
\label{sigma>1.sec}

When $\sigma \geq 1$, each machine only need to output a one-bit measurement to achieve the global optimal rate as if there are no communication constraints. Some related results are available in \cite{kipnis2019mean}. The following procedure is based on the setting when $b_i=1$ for all $i=1,...,m$. If $b_i > 1$ for some $i$, then one can simply discard all remain bits so that only one bit is sent by each machine. 

Here is the MODGAME procedure when $\sigma \geq 1$:
\begin{enumerate}
    \item[Step 1.] The $i-$th machine outputs
    \[ Z_i = \begin{cases} 0 \quad \text{if } X_i < 0\\ 1 \quad \text{if } X_i \geq 0 \end{cases}.    \]
    
    \item[Step 2.] The central machine collects $Z_1, Z_2, ..., Z_m$ and estimates $\theta$ by
    \[  \hat\theta_D = \tau_{[0,1]}\left(\sigma\Phi^{-1}\left(\frac{1}{m}\sum_{i=1}^m Z_i \right)\right)   \]
    
    where $\tau$ is the truncation function defined in \eqref{def_tau} and $\Phi$ is the cumulative distribution function of a standard normal,
    $\Phi(x) \triangleq \frac{1}{\sqrt{2\pi}}\int_{-\infty}^{x} e^{t^2/2}dt$. Here $\Phi^{-1}$ is the inverse of $\Phi$ and we extend it by defining $\Phi^{-1}(0) = -\infty$ and $\Phi^{-1}(1) = \infty$.
\end{enumerate}

\subsection{Theoretical properties of the MODGAME procedure}
\label{upperbd.sec}

Section \ref{MODGAME.sec} gives a detailed construction of the MODGAME procedure, which clearly satisfies the communication constraints by construction. The following result provides a theoretical guarantee for the statistical performance of MODGAME.

\begin{theorem}
	\label{unified_upperbd.thm}
    For given communication budgets $b_{1:m}$ with $b_i \ge 1$ for $i=1, \ldots, m$, let $B = \sum_{i=1}^m b_i$ and let $\hat\theta_D$ be the  MODGAME estimate. Then there exists a constant $C>0$ such that
	\beq \label{eq_ub} \sup_{\theta\in [0,1]} \E(\hat\theta_D - \theta)^2 \leq \begin{cases} C \cdot 2^{-2B} &\text{ if } B<\log\frac{1}{\sigma} + 2 \\ C \cdot \frac{\sigma^2}{(B-\log\frac{1}{\sigma})} &\text{ if } \log\frac{1}{\sigma} + 2 \leq B < \log\frac{1}{\sigma} + m \\ C \cdot \left(\frac{\sigma^2}{m} \wedge 1\right) &\text{ if } B \geq \log\frac{1}{\sigma} + m  \end{cases} . \eeq
\end{theorem}
An interesting and somewhat surprising feature of the upper bound is that it depends on the communication constraints $b_{1:m}$ only through the total budget $B = \sum_{i=1}^m b_i$, not the specific value of $b_{1:m}$, so long as each machine can transmit at least one bit.

\subsection{Lower bound analysis and discussions}
\label{lowerbd.sec}

Section \ref{MODGAME.sec} gives a detailed construction of the MODGAME procedure and Theorem \ref{unified_upperbd.thm} provides a theoretical guarantee for   the estimator. 
We shall now prove that  MODGAME is indeed rate optimal among all estimators satisfying the communication constraints by showing that the upper bound in Equation \eqref{eq_ub} cannot be improved. More specifically, the following lower bound provides a fundamental limit on the estimation accuracy under the communication constraints.

\begin{theorem}
	Suppose $b_i \geq 1$ for all $i = 1,2,...,m$. Let $B = \sum_{i=1}^m b_i$. Then there exists a constant $c>0$ such that
	\[ R_1(b_{1:m})  \geq \begin{cases} c \cdot 2^{-2B} &\text{ if } B< \log \frac{1}{\sigma} + 2 \\ c \cdot \frac{\sigma^2}{(B-\log \frac{1}{\sigma})} &\text{ if } \log \frac{1}{\sigma} + 2 \leq B < \log \frac{1}{\sigma} + m \\ c \cdot \left(\frac{\sigma^2}{m} \wedge 1\right) &\text{ if } B \geq \log \frac{1}{\sigma} + m. \end{cases} \]
	\label{lower-bound}
\end{theorem}

The key novelty in the lower bound analysis is the decomposition of the statistical risk into \textit{localization} error and \textit{refinement} error. By assigning $\theta$ a uniform prior on the candidate set
\[ G_\delta \triangleq \left\{ \theta_{u,v} = \sigma u + \delta v: u = 0,1,2,...,\left(\lfloor \frac{1}{\sigma} \rfloor - 1\right), v = 0,1  \right\}, \]
where $\delta$ is a precision parameter that will be specified later, estimation of $\theta$ can be decomposed into estimation of $u$ and $v$. One can view estimation of $u$ as the \textit{localization} step and estimation of $v$ as the \textit{refinement} step. The following lemma is a key technical tool.

\begin{lemma} \label{lm-decomposition}
	Let $0<\sigma<1$ and let  $u$ be uniformly distributed on $\{0,1,...,\lfloor \frac{1}{\sigma} \rfloor - 1\}$ and $v$ be uniformly distributed on $\{0,1\}$. Let $u$ and $v$ be independent and let $\theta = \theta_{u,v} = \frac{u}{\sigma} + \delta v$ where $0 < \delta < \frac{\sigma}{8}$. Then for all $\hat\theta \in \mathcal{A}(b_{1:m})$, 
	\beq I(\hat\theta; u) + \frac{\sigma^2}{64\delta^2}I(\hat\theta; v) \leq B.  \label{lm-decomposition-ineq} \eeq
\end{lemma}

\begin{remark} \rm
	The proof of Lemma \ref{lm-decomposition} mainly relies on the strong data processing inequality (Lemma 14 in \cite{CaiWei2019Supplement}).
	The strong data processing inequality was originally developed in information theory, for reference see \cite{raginsky2016strong}.  \cite{zhang2013information} and \cite{braverman2016communication} applied this technical tool to obtain lower bounds for distributed mean estimation. However, their lower bounds are not sharp in general, due to the fact that the focus was on bounding the refinement error using the strong data processing inequality, but failed to bound the localization error.
\end{remark}

Lemma \ref{lm-decomposition} suggests that under the communication constraints $b_{1:m}$, there is an unavoidable trade-off between the mutual information $I(\hat\theta; u)$ and $I(\hat\theta; v)$. So one of the above two quantities must be ``small". When $I(\hat\theta; u)$ (or $I(\hat\theta; v)$) is small than a certain threshold, it can be shown that the estimator $\hat\theta$ cannot accurately estimate $u$ (or $v$), which means the \textit{localization} error (or the \textit{refinement} error) is large. Given that one of \textit{localization} error and \textit{refinement} error must be larger than a certain value,  the desired lower bound follows. A detailed proof of Theorem \ref{lower-bound} is given in Section \ref{sec.proof}. 

\medskip\noindent
{\it Minimax rate of convergence.} Theorems \ref{unified_upperbd.thm} and \ref{lower-bound} together yield a sharp minimax rate for distributed univariate Gaussian mean estimation:
\beq 
\inf_{\hat\theta\in\mathcal{A}(b_{1:m})}\sup_{\theta\in [0,1]} \E(\hat\theta - \theta)^2 \asymp \begin{cases} 2^{-2B} &\text{ if } B< \log \frac{1}{\sigma} + 2 \\ \frac{\sigma^2}{(B-\log \frac{1}{\sigma})} &\text{ if } \log \frac{1}{\sigma} + 2 \leq B < \log \frac{1}{\sigma} + m \\ \frac{\sigma^2}{m} \wedge 1 &\text{ if } B \geq \log \frac{1}{\sigma} + m \end{cases}. 
\label{minimax_rate_univar} 
\eeq
The results also show that MODGAME is  rate optimal.

The minimax rate only depends on the total communication budgets $B = \sum_{i=1}^m b_i$. As long as each transcript contains at least one bit, how these communication budgets are allocated to local machines does not affect the minimax rate. This surprising phenomenon is due to the symmetry among the local machines since samples on different machines are independent and identically distributed. 


\begin{remark}
Figure \ref{pic_univ_rate} gives an illustration for the minimax rate $\eqref{minimax_rate_univar}$, which is divided into three phases: localization, refinement, and optimal-rate. The minimax risk decreases quickly in the localization phase, when the communication constraints are extremely severe; then it decreases slower in the refinement phase, when there are more communication budgets; finally the minimax rate coincides with the centralized optimal rate \citep{bickel1981minimax} and stays the same, when there are sufficient communication budgets. The value for each additional bit decreases as more bits are allowed. 
    
    In the localization phase, the risk is reduced to as small as $O(\sigma^2)$, which can be achieved by using the sample on only ONE machine and there is no need to ``communicate" with multiple machines. In the refinement phase, the risk is further reduced to $O(\sigma^2/m)$. However, one must aggregate information from all machines in order to achieve this rate.
\end{remark} 

\begin{remark}
	If the central machine itself also has an observation, or equivalently if one of the local machines serves as the central machine, then the communication constraints can be viewed as one of $b_i$ is equal to infinity. This setting is considered in some related literature, for instance, see \cite{jordan2018communication}. Then according to Theorem \ref{unified_upperbd.thm}, MODGAME always achieves the centralized rate $\frac{\sigma^2}{m}\wedge 1$, as long as at least one bit is allowed to communicate with each local machine. 
\end{remark}   
 
\begin{remark}
Our analysis on the minimax rate can be generalized to the $l_r$ loss for any $r \geq 1$, with suitable modifications on both the lower bound analysis and optimal procedure. 
\end{remark}

\section{Distributed Multivariate Gaussian Mean Estimation} \label{sec.multi}

We turn in this section to distributed estimation of a multivariate Gaussian mean under the communication constraints. Similar to the univariate case, suppose we observe on $m$ local machines  i.i.d. random samples:
\[  X_i \iid N_d(\theta,\sigma^2I_d), \quad i=1, \ldots, m, \]
where the $i$-th machine has access to $X_i$ only. Here we consider the multivariate Gaussian location family
\[  \mathcal{P}_\sigma^d = \left\{ N_d(\theta,\sigma^2I_d) : \theta \in [0,1]^d  \right\},  \]
where $\theta \in [0,1]^d$ is the mean vector of interest and the noise level $\sigma$ is known.
Under the constraints on the communication budgets $b_{1:m}$ with $b_i \ge 1$ for $i=1, \ldots, m$, the goal is to optimally estimate the mean vector $\theta$ under the squared error loss. 
We are interested in the minimax risk for the distributed protocol $\mathcal{A}(b_{1:m})$:
\[ R_d(b_{1:m}) =  \inf_{\hat\theta\in\mathcal{A}(b_{1:m})}\sup_{\theta\in [0,1]^d} \E \Vert \hat\theta - \theta \Vert^2. \]
Another goal is to develop a rate-optimal estimator that satisfies the communication constraints. The multivariate case is similar to the univariate setting, but it also has some distinct features, both in terms of the estimation procedure and the lower bound argument.
 
\subsection{Lower bound analysis}
\label{lowerbd-muivariate.sec}

We first obtain the minimax lower bound which is instrumental in establishing the optimal rate of convergence. The following lower bound on the minimax risk shows a fundamental limit on the estimation accuracy when there are communication constraints.  In view of the upper bound to be given in Section \ref{upperbd-muivariate.sec} that is achieved by a generalization of the  MODGAME procedure, the lower bound is rate optimal.

\begin{theorem} \label{multi-lb}
	Suppose $b_i \geq 1$ for all $i = 1,2,...,m$. Let $B = \sum_{i=1}^m b_i$ and $m' = \frac{1}{d}\sum_{i=1}^m (b_i\wedge d)$, then there exists a constant $c>0$ such that
	\[ R_d(b_{1:m})  \geq \begin{cases} c\cdot2^{-2B/d}d &\text{ if } B/d < \log \frac{1}{\sigma} + 2 \\ c\cdot\frac{d\sigma^2}{(B/d-\log \frac{1}{\sigma})} &\text{ if } \log\frac{1}{\sigma} + 2 \leq B/d < \log \frac{1}{\sigma} + (m'\vee 2) \\ c\cdot d\left(\frac{ \sigma^2}{m'}\wedge 1\right) &\text{ if } B/d \geq \log \frac{1}{\sigma} + (m'\vee 2) \end{cases}.  \]
\end{theorem}

Compared with the proof for the lower bound given in Theorem \ref{lower-bound} in the univariate setting, additional technical tools are needed to prove the lower bound $c\cdot d\left(\frac{ \sigma^2}{m'}\wedge 1\right)$ for the last case. A detailed proof of Theorem \ref{multi-lb} is given in the Supplementary Material \citep{CaiWei2019Supplement}.

\subsection{Optimal procedure}
\label{upperbd-muivariate.sec}

We  now construct an estimator of the mean vector under the communication constraints. Roughly speaking, the procedure, called multi-MODGAME,  first divides the communication budgets evenly into $d$ parts and then each part of communication budgets will be used to estimate one coordinate of $\theta$. Our analysis shows that multi-MODGAME achieves the minimax optimal rate under the communication constraints. The construction of the distributed estimator  $\hat\theta_D$ is divided into three steps.

\begin{figure}[htb]
	\centering
	\includegraphics[width=0.96\textwidth]{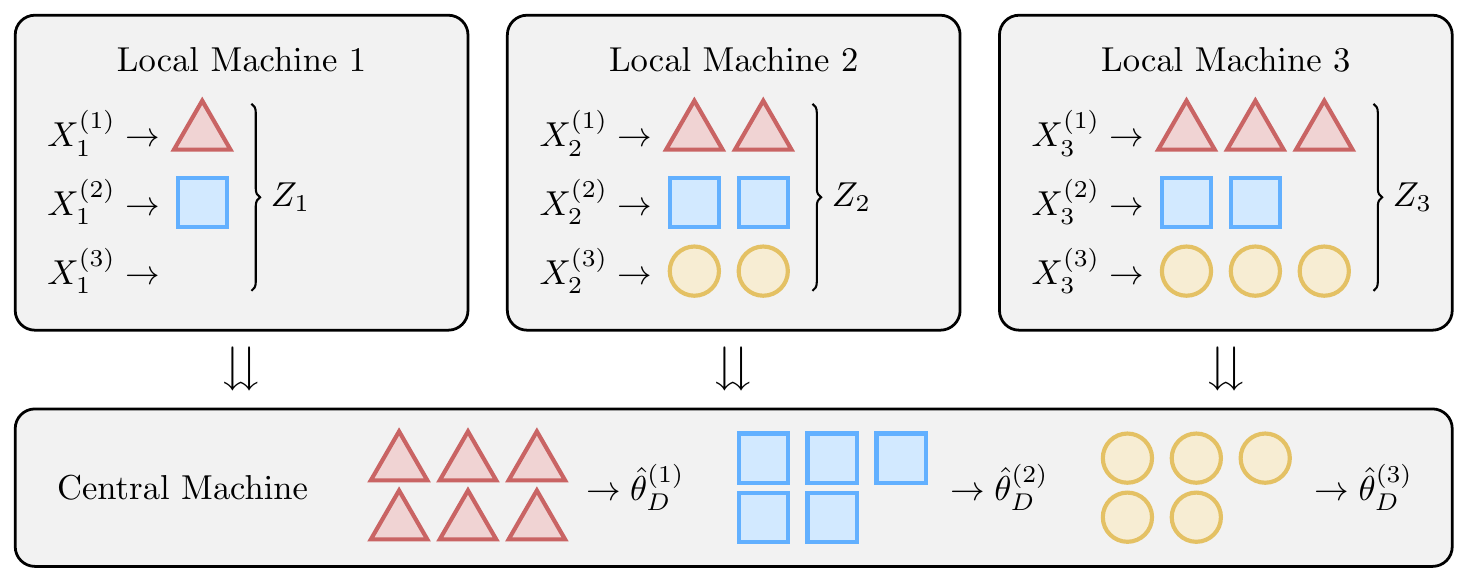}
	\caption{\small An illustration for multi-MODGAME. Communication budgets are evenly divided into three parts with each part used for estimating a coordinate of $\theta$ by the MODGAME procedure.  }
	\label{pic_multi_modgame}
	\vspace{-15pt}
\end{figure}

\medskip\noindent
\textbf{Step 1:} {\it Assign communication budgets.}
In this step we will calculate $\bp{k}_i \  (i=1,2,...,m; k=1,2,...,d)$ so that
\[  b_i = \bp{1}_i + \bp{2}_i + ... + \bp{d}_i \quad \text{for all } i=1,2,...m .  \]
where $\bp{k}_i$ is the number of bits within the transcript $Z_i$ which is associated with estimation of $\htp{k}$. 

Without loss of generality we assume $b_1 \leq b_2 \leq ... \leq b_m$, which can always be achieved by permuting the indices of the machines. 
Write $1,2,3,...,d$ repeatedly to form a sequence:
\[  Q \triangleq 1,2,3,...,d,1,2,3,...,d,1,2,3,...   \]

The sequence $Q$ is then divided into subsequences of lengths $b_1, b_2, ..., b_m$. Let $Q_1$ be the subsequence of $Q$ from index 1 to index $b_1$; let $Q_2$ be the next subsequence from index $b_1 + 1$ to $b_1 + b_2$; ... let $Q_m$ be the subsequence from index $\sum_{i=1}^{m-1} b_i + 1$ to $\sum_{i=1}^m b_i$.  For each $1\leq k \leq d$, let $\bp{k}_i$ be the number of occurrence of $k$ within $Q_i$. To be more precise, $\bp{k}_i$ can be calculated by
\[ \bp{k}_i = \left\lfloor \frac{\sum_{j=1}^i b_j - k}{d} \right\rfloor - \left\lfloor \frac{\sum_{j=1}^{i-1} b_j - k}{d} \right\rfloor  .  \]

\medskip\noindent
\textbf{Step 2:} {\it Generate transcripts on local machines.}
On the $i$-th machine, the transcript $Z_i$ is concatenated by short transcripts $Z_i^{(1)}, Z_i^{(2)}, ..., Z_i^{(d)}$, where the length of $Z_i^{(k)}$ is $\bp{k}_i$ for $k=1,2,...,d$.
Note that the $k$-th coordinate of the observations on each machine, $\Xp{k}_1, \Xp{k}_2, ..., \Xp{k}_m$, can be viewed as i.i.d univariate Gaussian variables with mean $\tp{k}$ and variance $\sigma^2$. For $1\leq k \leq d$, the transcripts $Z_1^{(k)}, Z_2^{(k)}, ..., Z_m^{(k)}$ can be generated the same way as if we implement MODGAME  to estimate $\tp{k}$ from observations $\Xp{k}_1, \Xp{k}_2, ..., \Xp{k}_m$, within the communication budgets $\bp{k}_1, \bp{k}_2, ..., \bp{k}_m$. Some machines may be assigned zero communication budget, if that happens those machines are ignored and  the procedure is implemented as if there are fewer machines.
 
\medskip\noindent
\textbf{Step 3:} {\it Construct distributed estimator $\hat\theta_D$.}
We have collected $Z_i^{(k)}$ ($i=1,2,...,m; k=1,2,...,d$) from the $m$ local machines. For $1\leq k \leq d$, as in MODGAME, one can use $Z_1^{(k)}, Z_2^{(k)}, ..., Z_m^{(k)}$ to obtain an estimate for $\htp{k}$:
\[  \hat\theta_D^{(k)} =  \hat\theta_D^{(k)}\left( Z_1^{(k)}, Z_2^{(k)}, ..., Z_m^{(k)}  \right) . \]
The final multi-MODGAME estimator $\hat\theta_D$ of the mean vector $\theta$  is just the vector consisting of the estimates for the $d$ coordinates:
\[  \hat\theta_D \triangleq \left(\hat\theta_D^{(1)}, \hat\theta_D^{(2)}, ..., \hat\theta_D^{(d)}  \right) . \]

The following result provides a theoretical guarantee for multi-MODGAME.

\begin{theorem}	\label{multi-ub}
	Let $B = \sum_{i=1}^m b_i$ and $m' = \frac{1}{d}\sum_{i=1}^m (b_i\wedge d)$. Then there exists a constant $C>0$ such that
	\beq  \sup_{\theta\in [0,1]^d} \E \Vert \hat\theta_D - \theta \Vert^2 \leq \begin{cases} C \cdot 2^{-2B/d}d &\text{ if } B/d < \log \frac{1}{\sigma} + 2 \\ C\cdot \frac{d\sigma^2}{(B/d-\log \frac{1}{\sigma})} &\text{ if } \log\frac{1}{\sigma} + 2 \leq B/d < \log \frac{1}{\sigma} + (m'\vee 2) \\ C\cdot d\left(\frac{ \sigma^2}{m'}\wedge 1\right) &\text{ if } B/d \geq \log \frac{1}{\sigma} + (m'\vee 2) . \end{cases} \label{multi_lb_eq} \eeq
\end{theorem}

The lower and upper bounds given Theorems \ref{multi-lb} and \ref{multi-ub} together establish the minimax rate for distributed multivariate Gaussian mean estimation:
\beq  \inf_{\hat\theta\in\mathcal{A}(b_{1:m})}\sup_{\theta\in [0,1]^d} \E \Vert \hat\theta - \theta \Vert^2 \asymp \begin{cases} 2^{-2B/d}d &\text{ if } B/d < \log \frac{1}{\sigma} + 2 \\ \frac{d\sigma^2}{(B/d-\log \frac{1}{\sigma})} &\text{ if } \log\frac{1}{\sigma} + 2 \leq B/d < \log \frac{1}{\sigma} + (m'\vee 2) \\ d\left(\frac{ \sigma^2}{m'}\wedge 1\right) &\text{ if } B/d \geq \log \frac{1}{\sigma} + (m'\vee 2) \end{cases} \label{minimax_rate_multi} \eeq
where $B = \sum_{i=1}^m b_i$ is the total communication budget and $m' = \frac{1}{d}\sum_{i=1}^m (b_i\wedge d)$ is the ``effective sample size".
In particular, the minimax rate \eqref{minimax_rate_univar} for the univariate case is an special case for the above minimax rate \eqref{minimax_rate_multi} with $d=1$.

\begin{remark}
Different from the univariate case, in the multivariate case the minimax rate depends on not only the total communication budget $B$, but also the effective sample size $m'$. How the communication budgets assigned to individual local machines affects the difficulty of the estimation problem.
If the communication budgets are tight on some machines, then one may have $m' \ll m$, which means the centralized minimax rate cannot be achieved even if the total communication budget $B$ is sufficiently large. 
\end{remark}
    
\begin{remark}
The present paper focuses on the unit hypercube $[0,1]^d$ as the parameter space.  A similar analysis can be applied to other ``regular" shape constraints, such as a ball or a simplex, and the minimax rate depends on the constraint. 
\end{remark}

\section{Simulation Studies} 
\label{sec.simulation} 
 
It is clear by construction that MODGAME and multi-MODGAME satisfy the communication constraints and are easy to implement. We investigate in this section their numerical performance through simulation studies. Comparisons with the existing methods are given and the results are consistent with the theory. 

We first consider MODGAME for estimating a univariate Gaussian mean. In this case, we set $d=1$ and $b_1=b_2=...=b_m=b$, i.e. the communication budgets for all machines are equal, and compare the empirical MSEs of MODGAME, naive quantization (see e.g. \cite{zhang2013information}), and  sample mean. For naive quantization, each machine projects its observation to $[0,1]$ and quantizes it to precision $2^{-b}$. The quantized observation is sent to the central machine and the central machine uses their average as the final estimate.
The  sample mean is the efficient estimate when there are no communication constraints, which can be viewed as a benchmark for any distributed Gaussian mean estimation procedure.

\begin{figure}[htb]
    \centering
    \begin{subfigure}[htb]{0.32\textwidth}
    	\includegraphics[width=\textwidth]{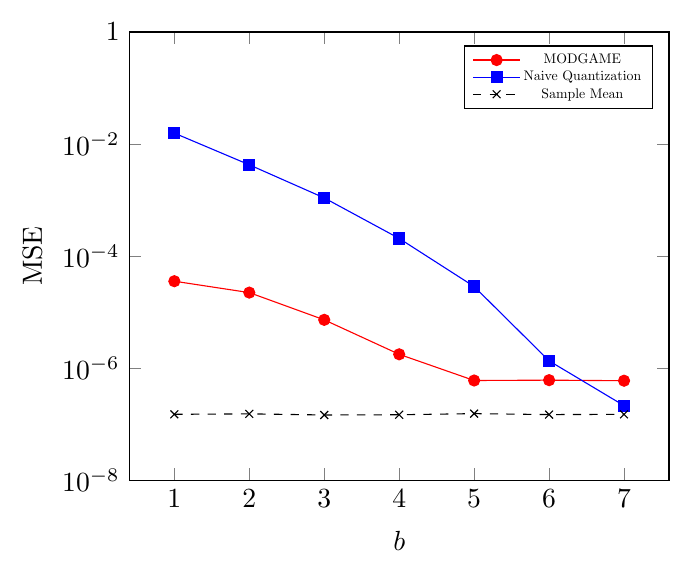}
    	\caption{\small Fixed $m$ and $\sigma$}
    	\label{fig:simulation_1}
    \end{subfigure}
    \begin{subfigure}[htb]{0.32\textwidth}
    	\includegraphics[width=\textwidth]{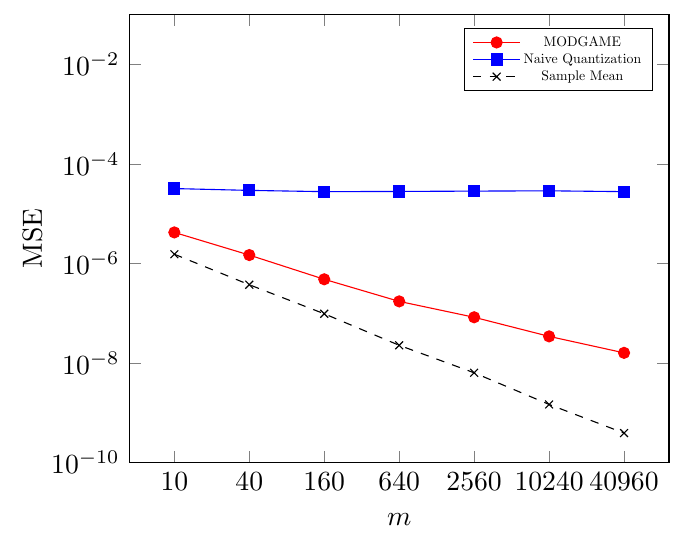}
    	\caption{\small Fixed $b$ and $\sigma$}
    	\label{fig:simulation_2}
    \end{subfigure}
	\begin{subfigure}[htb]{0.32\textwidth}
		\includegraphics[width=\textwidth]{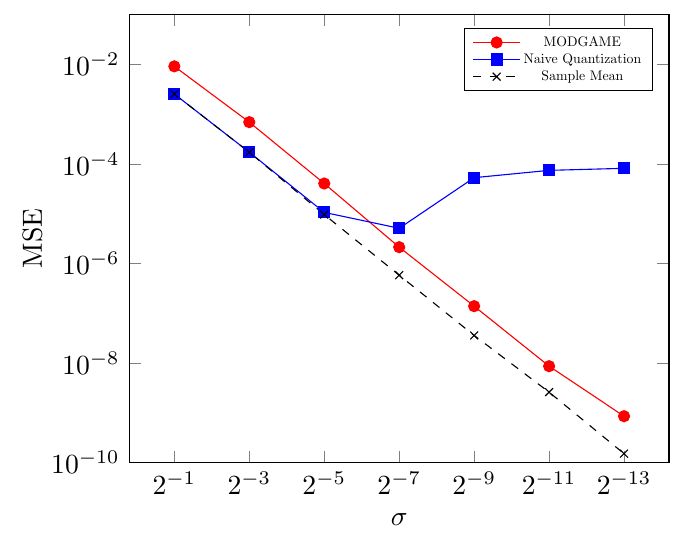}
		\caption{\small Fixed $b$ and $m$}
		\label{fig:simulation_3}
	\end{subfigure}
    \caption{\small Comparisons of the MSEs of MODGAME (red), naive quantization (blue) and sample mean (black). MSEs are plotted on log-scale. In \ref{fig:simulation_2} and \ref{fig:simulation_3}, $m$ and $\sigma$ are plotted on log-scale.} \vspace{-12pt}
\end{figure}

First, we fix $m=100$, $\sigma=2^{-8}$ and assign the communication budget for each machine $b$ from 1 to 7. 
The MSEs of  the three estimators are shown in Figure \ref{fig:simulation_1}, which shows that MODGAME makes better use of the communication resources in comparison to naive quantization.It can be seen from the figure, MODGAME outperforms naive quantization when the communication constraints are extremely severe. As the communication budgets increases, naive quantization can nearly achieve the optimal MSE, meanwhile MODGAME still performs very well. 


In the second setting, we fix $\sigma=2^{-8}$, $b=5$ and vary the number of machines $m$ from $10$ to $40960$.  Figure \ref{fig:simulation_2} plots the MSEs of the three methods. The MSE of MODGAME decreases as number of machine increases and outperforms naive quantization; the MSE of naive quantization remains constant as the quantization error plays a dominant role in the MSE.

Finally, we fix $b=5$, $m=100$ and vary the standard deviation $\sigma$ from $2^{-1}$ to $2^{-13}$. Figure \ref{fig:simulation_3} shows the MSEs of the three estimators. It can be seen that MODGAME is robust for all choices of $\sigma$. The difference between the MSE of MODGAME and the optimal MSE for non-distributed sample mean is small.  For naive quantization, it is as good as the optimal non-distributed sample mean when $\sigma$ is large. However, as seen in the previous experiment, when $\sigma$ is small, the MSE of naive quantization is dominated by the quantization error and is much larger than the MSE of MODGAME. 
In all three settings, it can be seen clearly that the MSE of  MODGAME decreases as the communication budgets increases. This is consistent with the theoretical results established in Section \ref{sec.univ}  and demonstrates the tradeoff between the communication costs and statistical accuracy.

We now turn to multi-MODGAME. Different values of the dimension $d$ yield similar phenomena. We use $d=50$ here for illustration.
When $d$ is larger than the number of bits that is allowed to communicate on each machine, naive quantization is not valid as  it is unclear how to quantize the  $d$ coordinates of the observed vector. As a comparison, it can be seen in the following experiments that multi-MODGAME still performs well even if  $d$ is large and the communication budgets are tight. 

\begin{figure}[htb]
	\centering
	\includegraphics[width=3in]{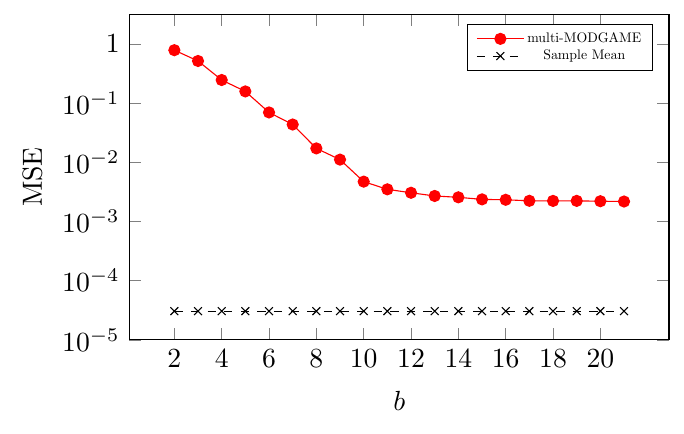}
	\caption{\small Comparisons of the MSEs of multi-MODGAME (red) and sample mean (black) where each machine is assigned $b=2,3,...,21$ bits. MSEs are plotted on log-scale.} \vspace{-12pt}
	\label{fig:simulation_4}
\end{figure}

Same as before, we set $b_1=b_2=...=b_m=b$, i.e. the communication budgets for all machines are equal. We set $d=50$, $\sigma=2^{-8}, m=25$ and assign the communication budgets $b$ for each machine from 2 to 21. The MSEs of multi-MODGAME  and sample mean are shown in Figure \ref{fig:simulation_4}. A phase transition at $b=10$ can be clearly seen. When $b\leq 10$, the MSE decreases quickly at an exponential rate. When $b>10$, the decrease becomes relatively slow. This phenomenon is consistent with the theoretical prediction that different phases appear in the convergence rate for multi-MODGAME (Theorem \ref{multi-ub}). 

\section{Discussion} 
\label{sec.discussion}

In this paper, we established the minimax rates of convergence for distributed Gaussian mean estimation under the communication constraints and developed two rate optimal estimation procedures, MODGAME for the univariate case and multi-MODGAME for the multivariate case. A key to solving this problem is the decomposition of the minimax estimation  problem into two steps, {\it localization} and {\it refinement}, which appears in both the lower bound analysis and optimal procedure design. 
 
In spite of these optimality results, there are still several open problems on distributed Gaussian mean estimation. For example, an interesting problem is the optimal estimation of the mean $\theta$ when the variance $\sigma^2$ is unknown. The lack of knowledge of  $\sigma^2$ requires additional communication efforts for optimally estimating $\theta$.  
When there are more than one sample available on each local machine, it is possible to estimate $\sigma^2$ locally on each machine and then use MODGAME (with some suitable modification) to estimate $\theta$. Another possibility is naive quantization introduced in Section \ref{sec.simulation}. 

Other than estimating the mean $\theta$, distributed estimation of the variance $\sigma^2$ is also an interesting and important problem. When there are multiple samples on each local machine, the local estimate of $\sigma^2$ can be viewed as an observation drawn from a $\chi^2$ distribution. The problem  then becomes a distributed $\chi^2$ estimation problem and it might be solved by using a similar approach to the one used in the present paper. We leave it for future work. 

Optimal estimation of the mean of a multivariate Gaussian distribution with a general (known) covariance matrix is another interesting problem. A naive approach is to ignore the dependency and apply MODGAME to estimate the coordinates individually, this is arguably not communication efficient in general.
For instance, if the correlation between certain coordinates is large, it may be possible to save a significant amount of communication budget by utilizing the information from one coordinate to help estimate the other. Another approach is to use multi-MODGAME after orthogonalization. More specifically,  consider the Gaussian location family with a general non-singular covariance matrix $\Sigma$. 
Let $\lambda_{\min}>0$ be the smallest eigenvalue of $\Sigma$. For $X \sim N_d(\theta, \Sigma)$, 
$\lambda_{\min}^{1/2}(d\Sigma)^{-1/2}X \sim  N_d\left(\lambda_{\min}^{1/2}(d\Sigma)^{-1/2}\theta, \frac{\lambda_{\min}}{d}I_d\right).$ 
Note that $\lambda_{\min}^{1/2}(d\Sigma)^{-1/2}\theta \in [0,1]^d$ for any $\theta \in [0,1]^d$, therefore one can apply multi-MODGAME to estimate $\lambda_{\min}^{1/2}(d\Sigma)^{-1/2}\theta$, then transform it back to get an estimate for $\theta$. However, this estimate is generally not rate-optimal. A systematic study is needed for this problem.
       
This paper arguably considered one of the simplest settings for optimal distributed estimation under the communication constraints, but as can be seen in the paper, both the  construction of the rate optimal estimators and the theoretical analysis are already quite involved for such a seemingly simple problem. As we deepen our understanding on distributed learning under the communication constraints, we hope to extent this line of work to investigate other statistical problems in distributed settings, including high-dimensional linear regression, large-scale multiple testing and nonparametric function estimation. 
In some of these problems, feedbacks (communications from the central machine to the local machines) appear to be necessary. It is interesting to understand fully when and to what extend feedbacks help in terms of improving statistical accuracy.

\section{Proofs} \label{sec.proof}

We prove two main results, Theorems  \ref{unified_upperbd.thm} and  \ref{lower-bound}, for the univariate case. For reasons of space, Theorems \ref{multi-lb} and \ref{multi-ub} and the technical lemmas are proved in the Supplementary Material \citep{CaiWei2019Supplement}.

\subsection{Proof of Theorem \ref{unified_upperbd.thm}}

We divide into two cases: $\sigma < 1$ and $\sigma \geq 1$.

\subsubsection{Proof of Theorem \ref{unified_upperbd.thm} when $\sigma < 1$}

We first define the ``change points sets" for the Gray functions $g_k(x)$ and conjugate Gray functions $\bar{g}_k(x)$. For any $k\geq 1$, let $G_k$ be the change points set for $g_k$, which is defined as
\[ G_k \triangleq \{ (2j-1) \cdot 2^{-k} : 1\leq j \leq 2^{k-1} \}.  \]

Similarly, let $\bar{G}_k$ be the change-points set for $\bar{g}_k$, which is defined as
\[ \bar{G}_k \triangleq \{ j \cdot 2^{k-1} : 1\leq j \leq 2^{k-1} - 1  \}.    \]

As the name suggests, the change points set for a Gray function (or a conjugate Gray function) is the collection of points $x\in [0,1] $ where $g_k(x)$ (or $\bar{g}_k(x)$) changes its value from 0 to 1 or from 1 to 0. More precisely, 
\[  G_k = \{x: \lim_{y\to x^-}g_k(y) \neq \lim_{y\to x^+}g_k(y) \} \quad {\rm and } \quad \bar{G}_k = \{x: \lim_{y\to x^-}\bar{g}_k(y) \neq \lim_{y\to x^+}\bar{g}_k(y) \}. \]

An important property for the change-points sets is that for any $k\geq 1$,
\beq  \bar{G}_{k+1} = \bigcup_{i=1}^{k} G_k \quad {\rm and} \quad G_i\cap G_j = \varnothing \  \forall 1\leq i < j \leq k. \label{prop_lattice} \eeq

\noindent
\textbf{Case 1: $B < \log\frac{1}{\sigma} + 2$. }
We first state several technical lemmas in general forms. These lemmas  will also be used in Case 2.

\begin{lemma} \label{lm_ub_1}
Let $x\in [0,1]$ and let $K \geq 1$ be a positive integer. Let $g_1, ..., g_K$ be the Gray functions and let $G_1, ..., G_K$ be the corresponding sets of change points.
Let $z_k = g_k(x)$, $k \in \{1,...,K\}$. If $d(x,G_k) > 2^{-(K+2)}$ for $k \in \{1,...,K\}$, then 
	\[ d(x, \Dec_K(z_1,z_2,...,z_K)) \leq 2^{-(K+2)}.   \]
\end{lemma}

\begin{lemma} \label{lm_ub_2}
	If $z_k = g_k(X)$ where $X \sim N(x,\sigma^2)$, then
	\beq \P(z_k \neq g_k(x)) \leq 2e^{-\frac{d(x, G_k)^2}{2\sigma^2}}. \label{eq_ub_lm2_1} \eeq
	
	Similarly, if $\bar{z}_k = g_k(X)$ where $X \sim N(x,\sigma^2)$, then
	\beq \P(\bar{z}_k \neq \bar{g}_k(x)) \leq 2e^{-\frac{d(x, \bar{G}_k)^2}{2\sigma^2}}.   \label{eq_ub_lm2_2} \eeq
\end{lemma}

\begin{lemma} \label{lm_ub_3}
	Fix any $x\in [0,1]$ and integer $1\leq K \leq \log\frac{1}{\sigma} + 2$. For any $1\leq k \leq K$, let $z_k = g_k(X_k)$ where $X_k \sim N(x,\sigma^2)$. ($X_1$,$X_2$,...,$X_K$ can be correlated.) Then there exists a constant $C_1>0$ such that, for any $L\leq K$, 
	\[  \P(d(x, \Dec_K(z_1,z_2,...,z_K))\geq \frac{5}{4}2^{-L} - 2^{-K}) \leq C_1e^{-\frac{2^{-2(L+2)}}{2\sigma^2}}.    \]
\end{lemma}

Now we prove Case 1. For simplicity denote $A = d(\theta, \Dec_B(U_1,U_2,...,U_B))$. Note that $A \leq 1$, so we have
\begin{align*}
	\E A^2 &\leq \P(A \leq \frac{5}{4}2^{-B})\cdot(\frac{5}{4}2^{-B})^2 + \sum_{k=0}^{B-1} \P(\frac{5}{4}2^{-B+k} \leq A \leq \frac{5}{4}2^{-B+k+1})\cdot(\frac{5}{4}2^{-B+k+1})^2 \\
	&\leq 1\cdot(\frac{5}{4}2^{-B})^2 + \sum_{k=0}^{B-1} \P(A \geq \frac{5}{4}2^{-B+k})\cdot(\frac{5}{4}2^{-B+k+1})^2.
\end{align*}

Note that $B \leq \log\frac{1}{\sigma} + 2$, and $U_k$ has the same distribution as $g_k(X)$ where $X \sim N(\theta,\sigma^2)$. We can apply Lemma \ref{lm_ub_3} and further get
\begin{align*}
	\E A^2 &\leq \frac{25}{16}2^{-2B} + \sum_{k=0}^{B-1} C_1e^{-\frac{2^{-2(B-k+2)}}{2\sigma^2}}\cdot\left(\frac{5}{4}2^{-B+k+1}\right)^2 \\
	&\leq \frac{25}{16}2^{-2B}\left( 1 + C_1 \sum_{k=0}^{B-1} 2^{-(2k+2)}e^{-\frac{2^{-2(-\log\sigma+2-k+2)}}{2\sigma^2}}  \right) \\
	&\leq C_2\cdot 2^{-2B},
\end{align*}
where $C_2 \triangleq \frac{25}{16}\left( 1 + C_1 \sum_{k=0}^{\infty} 2^{-(2k+2)}e^{-2^{(2k-9)}} \right)$ is summable.

Finally, we have $\hat\theta_D \in \Dec_B(U_1,U_2,...,U_B)$ and note that the length of $\Dec_B(U_1,U_2,...,U_B)$ is $2^{-B}$, therefore we conclude that
\[
	\E(\hat\theta_D - \theta)^2 \leq \E(A + 2^{-B})^2 \leq 2\E A^2 + 2^{-2B} \leq (2C_2+1)2^{-2B}.
\]
The upper bound in \eqref{eq_ub} for Case 1 is proved.

\noindent
\textbf{Case 2:} $\log\frac{1}{\sigma} + 2 \leq B < \log\frac{1}{\sigma} + m$.
We define
\beq \tilde{I}_1 = \{x: d(x, I_1') \leq 2^{-(\floorsig - \floorn) } \label{def_tI1} \} \cap [0,1]. \eeq
which is the interval that stretches out $\frac{1}{4}$ the length of $I_1'$ on both sides.

The proof is divided into three steps with each step summarized as a lemma below. These lemmas also imply the purpose of constructing intervals $I_1$ and $I_2$: they are confidence intervals with small risks of $\theta$ falling outside.  
\begin{lemma} \label{lm_ub_4}
	There exists a constant $C_3>0$ such that
	\[  \E((\hat\theta_D - \theta)^2\ind{\theta\notin \tilde{I}_1}) \leq \frac{C_3\sigma^2}{n}.   \]
\end{lemma}

\begin{lemma} \label{lm_ub_5}
	The set $I_2'$ defined in \eqref{def_I2p} is an interval and there exists a constant $C_4>0$ such that
	\[  \E((\hat\theta_D - \theta)^2\ind{\theta \in \tilde{I}_1, \theta\notin I_2}) \leq \frac{C_4\sigma^2}{n}.   \]
\end{lemma}

\begin{lemma} \label{lm_ub_6}
	(1) One of the following two conditions must hold:
	\[ I_2 \subseteq \left[(2j-\frac{3}{4})\cdot2^{-(\floorsig - 6)}, (2j+\frac{3}{4})\cdot2^{-(\floorsig - 6)}\right] \text{ for some } j\in \integers \] 
	or 
	\[ I_2 \subseteq \left[(2j+\frac{1}{4})\cdot2^{-(\floorsig - 6)}, (2j+\frac{7}{4})\cdot2^{-(\floorsig - 6)}\right] \text{ for some } j\in \integers . \]
	
	(2) There exists a constant $C_5>0$ such that
	\[  \E((\hat\theta_D - \theta)^2\ind{\theta \in I_2}) \leq \frac{C_5\sigma^2}{n}.   \]
\end{lemma}

From the above three lemmas we get
\begin{align*}
	\E((\hat\theta_D - \theta)^2 &\leq \E((\hat\theta_D - \theta)^2\ind{\theta\notin \tilde{I}_1}) + \E((\hat\theta_D - \theta)^2\ind{\theta \in \tilde{I}_1, \theta\notin I_2}) + \E((\hat\theta_D - \theta)^2\ind{\theta \in I_2}) \\
	&\leq (C_3+C_4+C_5)\frac{\sigma^2}{n}.
\end{align*} 

By the definition of $n$ in \eqref{def_n}, and $n\geq 1$, we know
\[
	B - \floorsig < \lfloor \log(n+1) \rfloor^2 + 2(n+1) < 2n + 2(n+1) \leq 6n.
\]   

Hence
\[ \E((\hat\theta_D - \theta)^2 \leq 6(C_3+C_4+C_5)\frac{\sigma^2}{B - \floorsig}  \leq 6(C_3+C_4+C_5)\frac{\sigma^2}{B - \log\frac{1}{\sigma}}.  \]

\noindent
\textbf{Case 3:} $B > \log\frac{1}{\sigma} + m$.
We can apply the procedure described in Case 2 (or Case 1 if $m=1$) as if we have $B' = \floorsig + m$ total communication budgets. So for some constant $C>0$ we have the guaranteed upper bound
\[ \E((\hat\theta_D - \theta)^2 \leq C\frac{\sigma^2}{B' - \log\frac{1}{\sigma}} \leq 2C\cdot \frac{\sigma^2}{m} \quad \text{if } m \geq 2 \]
or $ \E((\hat\theta_D - \theta)^2 \leq C\cdot 2^{-2B'} \leq C\cdot\frac{\sigma^2}{m}$  if $ m = 1$.  \qed

\subsubsection{Proof of Theorem \ref{unified_upperbd.thm} when $\sigma \geq 1$}

When $\sigma \geq 1$, we have $B \geq m \geq \log\frac{1}{\sigma} + m$, thus we only need to prove the last case in \eqref{eq_ub}, i.e.
\beq \sup_{\theta\in [0,1]} \E(\hat\theta_D - \theta)^2 \leq C \cdot \left(\frac{\sigma^2}{m} \wedge 1\right). \label{eq_ub_large_sigma} \eeq

Note that $\hat\theta_D, \theta \in [0,1]$ thus 
\beq  \E(\hat\theta_D - \theta)^2 < 1. \label{eq_ubls_1}  \eeq

Since $Z_1,Z_2,...,Z_m$ are $i.i.d$ Bernoulli with mean $P(X_i > 0) = \Phi(\frac{\theta}{\sigma})$,
\[  \E\left(\frac{1}{m}\sum_{i=1}^m Z_i - \Phi(\frac{\theta}{\sigma})\right)^2  \leq \frac{1}{4m}.   \]

Let $f(x) = \tau_{[0,1]}\left(\sigma\Phi^{-1}(x)\right)$. It is easy to verify that $f(x)$ is a $C\sigma$-Lipschitz function with some constant $C>0$ when $\sigma > 1$. Therefore we have
\[
	\E(\hat\theta_D - \theta)^2 
	\leq C^2\sigma^2 \E\left(\frac{1}{m}\sum_{i=1}^m Z_i - \Phi(\frac{\theta}{\sigma})\right)^2 
	\leq \frac{C^2}{4}\cdot\frac{\sigma^2}{m}.
\]
Inequality \eqref{eq_ub_large_sigma} follows by combining the above inequality and \eqref{eq_ubls_1}. \qed

\subsection{Proof of Theorem \ref{lower-bound}}

The lower bound is established separately for the  three cases: $B < \log \frac{1}{\sigma} + 2$, $\log \frac{1}{\sigma} + 2 \leq B < \log \frac{1}{\sigma} + m$, and $B \geq \log \frac{1}{\sigma} + m$. We shall first focus on the most important case $\log \frac{1}{\sigma} + 2 \leq B < \log \frac{1}{\sigma} + m$. The other two cases are relatively easy.  New technique tools are developed in the proof of this case.

\medskip\noindent
\textbf{Case 1:} $\log \frac{1}{\sigma} + 2 \leq B < \log \frac{1}{\sigma} + m$. Note that $b_i \geq 1$ for all $i = 1,2,...,m$ implies that $B = \sum_{i=1}^m b_i \geq m$. Therefore in this case we must have $\sigma < 1$.

Let $0 < \delta < \frac{1}{8}\sigma$ be a parameter to be specified later.  Define a grid of candidate values of $\theta$ as
\beq G_\delta \triangleq \left\{ \theta_{u,v} = \sigma u + \delta v: u = 0,1,2,...,\left(\lfloor \frac{1}{\sigma} \rfloor - 1\right), v = 0,1  \right\}. \label{def_grid} \eeq
Let $\U(G_\delta)$ be a uniform prior of $\theta$ on $G_\delta$. Note that $G_\delta \subset [0,1]$, so the minimax risk is lower bounded by the Bayesian risk:
\beq \inf_{\hat\theta\in\mathcal{A}(b_{1:m})}\sup_{\theta\in [0,1]}(\hat\theta - \theta)^2 \geq \inf_{\hat\theta\in\mathcal{A}(b_{1:m})}\E_{\theta\sim \U(G_\delta)}(\hat\theta - \theta)^2. \label{lb_c2_1}  \eeq

For any estimator $\hat\theta \in \mathcal{A}(b_{1:m})$, the rounded estimator $\hat\theta' \triangleq \argmin_{\tilde{\theta}\in G_\delta}|\tilde{\theta}-\hat\theta|$ always satisfy $(\hat\theta - \theta)^2 \geq \frac{1}{4}(\hat\theta'-\theta)^2$ for all $\theta\in G_\delta$. Note that $\hat\theta'$ also belongs to the protocol class $\mathcal{A}(b_{1:m})$, and only takes value in $G_\delta$, this implies
\beq \inf_{\hat\theta\in\mathcal{A}(b_{1:m})}\E_{\theta\sim \U(G_\delta)}(\hat\theta - \theta)^2 \geq \frac{1}{4}\inf_{\hat\theta\in\mathcal{A}(b_{1:m})\cap G_\delta}\E_{\theta\sim \U(G_\delta)}(\hat\theta - \theta)^2, \label{lb_c2_2} \eeq
where $\mathcal{A}(b_{1:m})\cap G_\delta$ is a shorthand for $\mathcal{A}(b_{1:m})\cap\{\hat\theta: \hat\theta \text{ only takes value in $G_\delta$}\}$.

Now we have $\hat\theta, \theta \in G_\delta$ thus they can be reparametrized by 
$\hat\theta = \theta_{\hat u, \hat v}$ and $ \theta = \theta_{u,v}$.
It is easy to verify the inequality
\[  (\hat\theta_{\hat u, \hat v} - \theta_{u,v})^2 \geq \max\left\{\frac{\sigma^2}{4}(\hat u - u)^2, \delta^2\ind{\hat v \neq v}  \right\}.  \]
Hence
\beq \inf_{\hat\theta\in\mathcal{A}(b_{1:m})\cap G_\delta}\E_{\theta\sim \U(G_\delta)}(\hat\theta - \theta)^2 \geq \inf_{\theta_{\hat u, \hat v}\in\mathcal{A}(b_{1:m})\cap G_\delta}\E_{\theta_{u,v}\sim \U(G_\delta)} \max\left\{\frac{\sigma^2}{4}(\hat u - u)^2, \delta^2\ind{\hat v \neq v} \right\}. \label{lb_c2_3} \eeq

Putting together \eqref{lb_c2_1}, \eqref{lb_c2_2}, and \eqref{lb_c2_3}, we have
\begin{equation}
\begin{aligned}
\inf_{\hat\theta\in\mathcal{A}(b_{1:m})}\sup_{\theta\in [0,1]}(\hat\theta - \theta)^2 \geq \frac{1}{4}\inf_{\theta_{\hat u, \hat v}\in\mathcal{A}(b_{1:m})\cap G_\delta}\E_{\theta_{u,v}\sim \U(G_\delta)} \max\left\{\frac{\sigma^2}{4}(\hat u - u)^2, \delta^2\ind{\hat v \neq v} \right\} 
\\ \geq \inf_{\theta_{\hat u, \hat v}\in\mathcal{A}(b_{1:m})\cap G_\delta} \max\left\{\frac{\sigma^2}{16}\E_{\theta_{u,v}\sim \U(G_\delta)}(\hat u - u)^2, \frac{\delta^2}{4}\P_{\theta_{u,v}\sim \U(G_\delta)}(\hat v \neq v) \right\}.
\end{aligned}
\label{lb_c2_4}
\end{equation}

Therefore, by assigning a prior $\theta\sim \U(G_\delta)$, we have successfully decomposed the estimation problem of $\theta$ into estimation problems of $u$ and $v$. We can view estimation of $u$ as ``localization" step and estimation of $v$ as ``refinement" step, so $\eqref{lb_c2_4}$ essentially has decomposed the statistical risk into localization error and refinement error.
To lower bound the right hand side of \eqref{lb_c2_4}, we show that under communication constraints, one cannot simultaneously estimate both $u$ and $v$ accurately, i.e. the localization and refinement errors cannot be both too small. Lemma \ref{lm-decomposition}, which shows that for any distributed estimator $\hat\theta$, there is unavoidable trade-off between the mutual information $I(\hat\theta; u)$ and $I(\hat\theta; v)$, is a key step.

We set $\delta = \frac{\sigma}{\sqrt{256(B+1-\log(\lfloor \frac{1}{\sigma} \rfloor))}}$, and assign the uniform prior $\U(G_\delta)$ to the parameter $\theta = \theta_{u,v}$. One can easily verify $\delta < \frac{1}{8}\sigma$, and $u,v$ are independent random variables where $u$ is uniform distributed on $\{0,1,...,\lfloor \frac{1}{\sigma} \rfloor - 1\}$, and $v$ is uniform distributed on $\{0,1\}$.
Therefore, we can apply Lemma \ref{lm-decomposition} to get inequality \eqref{lm-decomposition-ineq}. From the inequality \eqref{lm-decomposition-ineq} we can further get, for any $\hat\theta \in \mathcal{A}(b_{1:m}) \cap G_\delta$, one of the following two inequalities
\[ I(\hat\theta;u) \leq \log(\lfloor \frac{1}{\sigma} \rfloor) - 1 \quad {\rm or} \quad I(\hat\theta;v) \leq \frac{64\delta^2}{\sigma^2}\left(B+1-\log(\lfloor \frac{1}{\sigma} \rfloor)\right) \]
must hold. We show that either of the above bounds on the mutual information will result in a large statistical risk. 

\medskip\noindent
\textbf{Case 1.1:} $I(\hat\theta;u) \leq \log(\lfloor \frac{1}{\sigma} \rfloor) - 1$.
Note that $\hat u$ is a function on $\hat\theta$, thus by data processing inequality,
$I(\hat u; u) \leq I(\hat\theta; u) \leq \log(\lfloor \frac{1}{\sigma} \rfloor) - 1.$
Note that $u$ is uniform distributed on  $\{0,1,...,\lfloor \frac{1}{\sigma} \rfloor - 1\}$, thus $H(u) = \log(\lfloor \frac{1}{\sigma}\rfloor)$. We have
\beq H(u|\hat u) = H(u) - I(\hat u; u) \geq 1. \label{lb_c2_5} \eeq

The following lemma shows that large conditional entropy will result in large $L_2$ distance between two integer-value random variables. 
\begin{lemma} \label{l1fano}
	Suppose $A, D$ are two integer-value random variables. If $H(A|D)\geq \frac{1}{2}$, then there exist a constant $c_2>0$ such that 
	\[ \E (A-D)^2 \geq c_2.  \]
\end{lemma}

Given \eqref{lb_c2_5} and the fact that $\hat u$, $u$ are integer valued,  Lemma \ref{l1fano} yields
\beq  \E_{\theta_{u,v} \sim \U(G_\delta)}(\hat u - u)^2  \geq c_2.  \label{lb_c2_c11}  \eeq

\medskip\noindent
\textbf{Case 1.2:} $I(\hat\theta;v) \leq \frac{\delta^2}{c_1\sigma^2}(B+1-\log(\lfloor \frac{1}{\sigma} \rfloor))$.
By the strong data processing inequality, plug in $\delta = \frac{\sigma}{\sqrt{256(B+1-\log(\lfloor \frac{1}{\sigma} \rfloor))}}$ we have $I(\hat v; v) \leq I(\hat\theta; v) \leq \frac{1}{4}$, so
\[ H(v | \hat v) = H(v) - I(\hat v; v) \geq \frac{3}{4}. \]

It follows from Lemma \ref{l1fano} that  
\beq \P_{\theta_{u,v} \sim \U(G_\delta)}(\hat v \neq v) = \E_{\theta_{u,v} \sim \U(G_\delta)}(\hat v - v)^2 \geq c_2. \label{lb_c2_c12} \eeq

Combine \eqref{lb_c2_c11} for Case 1.1 and \eqref{lb_c2_c12} for Case 1.2 together, we have for any $\hat\theta \in \mathcal{A}(b_{1:m}) \cap G_\delta$,
\begin{equation}
\begin{aligned}
&\max\left\{\frac{\sigma^2}{16}\E_{\theta_{u,v}\sim \U(G_\delta)}(\hat u - u)^2, \frac{\delta^2}{4}\P_{\theta_{u,v}\sim \U(G_\delta)}(\hat v \neq v) \right\} \\ \geq& c_2\min\left\{ \frac{\sigma^2}{16}, \frac{\delta^2}{4} \right\}
= \frac{ c_2\sigma^2 }{1024(B+1-\log(\lfloor \frac{1}{\sigma} \rfloor))} \geq \frac{c_2}{2048}\cdot\frac{\sigma^2}{(B-\log\frac{1}{\sigma})}\ .
\end{aligned}
\label{lb_c2_6}
\end{equation}    

The minimax lower bound follows by combining \eqref{lb_c2_4} and \eqref{lb_c2_6}, 
\[ \inf_{\hat\theta\in\mathcal{A}(b_{1:m})}\sup_{\theta\in [0,1]}(\hat\theta - \theta)^2 \geq \frac{c_2}{2048}\cdot\frac{\sigma^2}{(B-\log\frac{1}{\sigma})} \ . \] 

\medskip\noindent
\textbf{Case 2:} $B < \log\frac{1}{\sigma} + 2$.
Let $S = 2^{B+1}$ and $K_S \triangleq \{\frac{i}{S}: i = 0,1,...,S-1\}$.
Denote by $\U(K_S)$ the uniform distribution on $K_S$. For the same reason as \eqref{lb_c2_1} and \eqref{lb_c2_2} we have
\begin{equation}
\begin{aligned}
\inf_{\hat\theta\in\mathcal{A}(b_{1:m})}\sup_{\theta\in [0,1]}(\hat\theta - \theta)^2 &\geq \inf_{\hat\theta\in\mathcal{A}(b_{1:m})}\E_{\theta\sim \U(K_S)}(\hat\theta - \theta)^2 \\ &\geq \frac{1}{4}\inf_{\hat\theta\in\mathcal{A}(b_{1:m})\cap K_S}\E_{\theta\sim \U(K_S)}(\hat\theta - \theta)^2 \\
&= \frac{1}{4S^2}\inf_{\hat\theta\in\mathcal{A}(b_{1:m})\cap K_S}\E_{\theta\sim \U(K_S)}(S\hat\theta - S\theta)^2 .
\end{aligned}
\label{lb_c1_1}
\end{equation} 

The parameter $\theta$ can be treated as a random variable drawn from the prior distribution $\U(K_S)$. Note that by the data processing inequality, for any $\hat\theta\in\mathcal{A}(b_{1:m})$,
\begin{align*}
	I(\hat\theta; \theta) &= I(\hat\theta(Z_1,Z_2,...,Z_m); \theta) 
	\leq I(Z_1,Z_2,...,Z_m; \theta) \\
	&\leq H(Z_1,Z_2,...,Z_m) 
	\leq \sum_{i=1}^m H(Z_i) \leq \sum_{i=1}^m b_i = B .
\end{align*} 

By $\theta \sim \U(K_S)$ we have
$ H(\theta | \hat\theta) = H(\theta) - I(\hat\theta; \theta) \geq \log S - B \geq 1$.
Note that when $\theta \sim \U(K_S)$, for any $\hat\theta \in \mathcal{A}(b_{1:m}) \cap K_S$, $S\hat\theta$ and $S\theta$ both take value in $\{0,1,2,...,S-1\}$. Also we have $H(S\theta | S\hat\theta) = H(\theta | \hat\theta) \geq 1$. Therefore, Lemma \ref{l1fano} yields that
$ \E_{\theta \sim \U(K_S)}(S\hat\theta - S\theta)^2 \geq c_2$.
We thus conclude that
\[ \frac{1}{4S^2}\inf_{\hat\theta\in\mathcal{A}(b_{1:m})\cap K_S}\E_{\theta\sim \U(K_S)}(S\hat\theta - S\theta)^2 \geq \frac{c_2}{4 \cdot 2^{2(B+1)}} = \frac{c_2}{16}\cdot 2^{-2B} .  \]
The desired lower bound follows by plugging into $\eqref{lb_c1_1}$.

\medskip\noindent
\textbf{Case 3:} $B \geq \log\frac{1}{\sigma} + m$.
The minimax risk for distributed protocols is always lower bounded by the minimax risk with no communication constraints:
\[  \inf_{\hat\theta\in\mathcal{A}(b_{1:m})}\sup_{\theta\in [0,1]}(\hat\theta - \theta)^2 \geq \inf_{\hat\theta}\sup_{\theta\in [0,1]}(\hat\theta - \theta)^2 \asymp \frac{\sigma^2}{m} \wedge 1 .  \]

which is given in \cite{bickel1981minimax}.  \qed

\begin{supplement}[id=supplemnt]
	\sname{Supplement A}
	\stitle{Supplement to ``\papertitle"}
	\slink[doi]{url to be specified}
	\sdescription{In this supplementary material, we provide detailed proofs for Theorems \ref{multi-lb} and \ref{multi-ub}, and proofs for technical lemmas.}
\end{supplement}

\bibliographystyle{apalike}
\bibliography{Distributed-Estimation-Reference.bib}


\end{document}